%% file: main.tex
\documentclass[11pt]{amsart}

\usepackage{amsmath}
\usepackage{amssymb}
\usepackage{graphicx}
\usepackage{bm}
\usepackage{graphicx}         %
\usepackage{caption}
\usepackage{subcaption}       %
\usepackage{rotating}         %

\newtheorem{theorem}{Theorem}[section]

\newtheorem{lemma}[theorem]{Lemma}
\newtheorem{proposition}[theorem]{Proposition}
\theoremstyle{definition}

\newtheorem{remark}[theorem]{Remark}

\numberwithin{equation}{section}

\usepackage{appendix}

\newcommand{\vect}[1]{\bm{#1}}     %

\usepackage{xcolor}                   %
\definecolor{mybackgroundblue}{RGB}{173, 216, 230} %
\definecolor{mytextblue}{RGB}{25, 25, 112}         %

\usepackage[color=yellow!50
, textcolor=mytextblue
, bordercolor=mytextblue
, colorinlistoftodos
,disable
]{todonotes}

\usepackage{enumerate}
\usepackage{enumitem}

\usepackage[hidelinks]{hyperref} %
\usepackage{cleveref}           %

\usepackage{varwidth}  %
\usepackage{lipsum}    %
\usepackage{pdfpages}  %
\usepackage{subfiles}

\title[Analysis of coupled BT models]{Analysis of linear Boussinesq-type models coupled with static interfaces}

\author[J. Galaz Mora]{José Galaz Mora}
\address[J. Galaz Mora]{CIGIDEN Centro de Investigación para la Gestión Integrada del Riesgo de Desastres, Pontificia Universidad Católica de Chile, Av. Vicuña Mackenna 4860, Santiago, Chile}
\email{{\tt jdgalazm@gmail.com}}

\author[M. Kazolea]{Maria Kazolea}
\address[M. Kazolea]{INRIA, Univ. Bordeaux, CNRS, Bordeaux INP, IMB, UMR 5251, 200 Avenue de la Vieille Tour, 33405 Talence cedex, France}
\email{\tt maria.kazolea@inria.fr}

\author[A. Rousseau]{Antoine Rousseau}
\address[A. Rousseau]{INRIA, Univ. Montpellier, CNRS, IMAG, UMR 5149, 860 Rue Saint Priest - 34095 Montpellier cedex 5, France}
\email{{\tt antoine.rousseau@inria.fr}}



\date{\today}

\begin{document}

\begin{abstract}
  We derive a new approach to analyze the coupling of linear Boussinesq and Saint-Venant shallow water wave equations in the case where the interface remains at a constant position in space. 

  We propose a one-way coupling model as a reference, which allows us to obtain an analytical solution, prove the well-posedness of the original coupled model and compute what we call the coupling error—a quantity that depends solely on the choice of transmission conditions at the interface. We prove that this coupling error is asymptotically small for a certain class of data and discuss its role as a proxy for the full error with respect to the 3D water wave problem. Additionally, we highlight that this error can be easily computed in other scenarios.
  
  We show that the coupling error consists of reflected waves and argue that this explains some previously unexplained spurious oscillations reported in the literature. Finally, we prove the well-posedness of the half-line linear Boussinesq problem.
\end{abstract}

\maketitle

\subfile{01-intro.tex}

\subfile{02-framework.tex}

\subfile{03-hybrid-analysis.tex}

\section{Conclusions}
\label{sec:conclusions}
We have studied the coupling of the linearized Boussinesq and Saint-Venant equations. By employing one-way models, we derived an analytical solution that allowed us to establish several key properties, including the well-posedness, a precise notion of the reflected waves, the quantification of their size and their link with the global error relative to the Euler equations, highlighting how they isolate the contribution of the chosen coupling conditions.

To the best of our knowledge, the approach used in this analysis is the first of its kind. Moreover, it can be readily extended to other linear equations, such as different BT models, their discretization and the implementation of other (better) coupling conditions, which we identify as the next steps in this research. 
\appendix
\section{Boussinesq half-line problem}
\subfile{A01halfline.tex}

\end{document}

%% file: 01-intro.tex
\section{Introduction}
\label{sec:ch5-linearmodels}

Shallow-water models such as Boussinesq-type  (BT) equations and Saint-Venant (nonlinear shallow water) equations are reduced models that provide asymptotic approximations of the full free surface Euler equations in shallow water. Their study has impacted important questions across disciplines including risk analysis, marine energy extraction, beach morphology evolution, etc. Many BT models exist,(see reviews \cite{brocchini2013reasoned, kazolea2024, coulaud2025comparison}) and much is known nowadays regarding their mathematical justification and properties (see ref. \cite{lannes2013water} for example). With this many available choices a path to build better models has been the coupling of BT equations, especially with the Saint-Venant (SV) equations. Such models can be better, on one side, in the sense of efficiency, by using expensive models only where necessary. This approach is common in heterogeneous domain decomposition  \cite{quarteroni1992heterogeneous,HeterogeneousDomain2008,gander2018multiscale,GanderMartin2023} and has already been done in tsunami modeling for example  \cite{IUGG1997,imamura2006tsunami,titov1997implementation,PedersenLøvholt2008}. On the other side, especially in the case of BT and SV coupling (BTSV), this coupling has allowed to represent wave breaking taking advantage of shock waves that naturally develop in nonlinear hyperbolic equations, despite having to neglect other important phenomena such as dispersive, nonhydrostatic, or others. Such an approach to coupling has become very popular due to its simplicity, efficiency and accuracy, as reviewed in \cite{kazolea2024}. More recently, it has been shown that one can obtain more robust models by coupling them conveniently, and benefit from properties that may be unknown for one model but well understood for another \cite{parisot2024thick}, such as absorbing and generating boundary conditions \cite{parisot2024thick,galaz2024coupling}, or better-behaved solutions in case of variable bathymetry or other data \cite{parisot2024thick}.

Coupled models have most commonly been justified by comparing them with the most expensive model everywhere and several techniques have been used for decades. See the short review in \cite{GanderMartin2023} and references therein, or the pioneering work in \cite{dubach1993contribution}. However, in our case, the reference model is the 3D Euler equations, hence the comparison can increase in complexity very easily and the model itself is not involved in the coupling. Moreover, several authors have observed multiple oscillations of different size and nature, sometimes just as inoquous reflections \cite{parisot2024thick} and \cite[Ch. 4]{galaz2024coupling}, even in models of blood flow \cite{quarteroni1992heterogeneous}, but other times as important issues that lead to instabilities and divergence of the numerical methods. See \cite{kazolea2018wave} and \cite[Ch. 1 and 4]{galaz2024coupling} for detailed reviews on these observations.  One important limitation of usual techniques presented for example in \cite{GanderMartin2023} is that they do not provide a way to understand the nature of the oscillations or a way to generalize it to more complex models. Regarding their well-posedness, only recently in \cite{parisot2024thick} the continuity of the model (and stability) on the initial data was proved with an energy method applied on a projected formulation of the BT models. To the authors' knowledge, so far there has been no discussion regarding how to measure their convergence to the 3D Euler equations and its relation with the oscillations, if any.

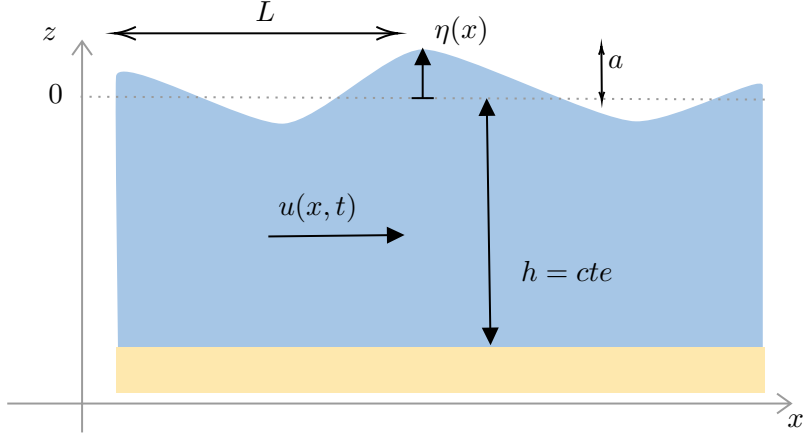
\begin{figure}  
    \centering 
    \input{figures/intro2/variables.tex}
    \caption{Sketch of variables and scales.}
    \label{fig:intro2-sketch}
\end{figure}

In this work we aim to contribute to the previous problem by studying the coupling of linearized Boussinesq and Saint-Venant equations with a static interface, \textit{i.e.} , we aim to find a function $\vect W(x,t) = (\eta(x,t), u(x,t))$  that  solves the linear Saint-Venant (SV) equations, which in adimensional form read
\begin{equation}
    \left\{
        \begin{array}{rl}
            \partial_t \eta +   \partial_x u &= 0 \\
            \partial_t u  +   \partial_x \eta &= 0
        \end{array}
        \right.
        \quad \text{ in } Q_{SV}
\end{equation}
and the linear Boussinesq (B) equations
\begin{equation}
    \left\{
        \begin{array}{rl}
            \partial_t \eta +   \partial_x u &= 0  \\
            \left(1 - \mu^2 \partial_x^2\right) \partial_t u  +  \partial_x \eta &= 0 \\ 
        \end{array}
        \right.   
        \quad \text{ in } Q_{B}
\end{equation}
with initial data
\begin{equation}
    \eta(\cdot,0) = \eta_0, \quad u(\cdot,0) = u_0,\quad  \text{ in }\Omega \
\end{equation}
where,  $\Omega = \mathbb{R}$, $T>0$ and  $Q = \Omega \times ]0,T[$. In the general case $Q_{SV}, Q_{B}$ are subdomains that form an arbitrary partition of $Q$, but here we consider the case of a static interface $x=0$, \emph{i.e.,} either $Q_{SV}=Q_-=\mathbb{R}^- \times]0,T[$ and $Q_{B}= Q_+ = \mathbb{R}^+ \times]0,T[$, or vice-versa. Also, as shown in \Cref{fig:intro2-sketch}, $\eta$ represents the free-surface elevation, $u$ the depth-averaged horizontal velocity, and $\mu>0$ is a parameter proportional to the ratio $h/L$ of (constant) depth $h>0$ to wavelength $L>0$. 

Letting
\begin{equation}
    A = \left(\begin{array}{cc}
        0 & 0     \\
        0 & \mu^2
    \end{array}\right),
    \quad
    B = \left(\begin{array}{cc}
        0 & 1 \\
        1 & 0
    \end{array}\right)
    \label{eq:ab-definition}
\end{equation}
one can write the problem more compactly as: Find $\vect W(x,t)$ such that

\begin{equation}
    \left\{
    \begin{array}{rl}
        \mathcal{L}_{SV}\vect W =0 & \quad \text{ in } Q_{SV} \\ 
        \mathcal{L}_B \vect W =0 & \quad \text{ in } Q_{B} \\
        \vect W(\cdot, 0) = \vect W_0 & \quad \text{ in } \Omega
    \end{array}\right.
    \label{eq:intro-abstract-coupling}
\end{equation}
where $\mathcal{L}_{B}=(I-A)\partial_t + B \partial_x$,  $\mathcal{L}_{SV}= \partial_t + B \partial_x$ and $\vect W_0=(\eta_0, u_0)$. 

Problem \cref{eq:intro-abstract-coupling} requires additional closure conditions to ensure a unique solution. Such closure could be achieved by decomposing the solution as 
\begin{equation}
    \vect W = \left\{
        \begin{array}{rl}
            \vect W_- & \text{ in } Q_- \\
            \vect W_+ & \text{ in } Q_+
        \end{array}
     \right.
\end{equation}
and
imposing boundary conditions with operators $B_\pm$ as
\begin{equation}
    \begin{split}
        B_-(\vect W_-, \vect W_+)=0 & \text{ at } \,  x=0\\ 
        B_+(\vect W_-, \vect W_+)=0 &  \text{ at } \, x=0        
    \end{split}
    \label{eq:intro-coupling-conditions}
\end{equation} but one could propose different ways on which the coupling could be done. One particular case of interest is the so-called hybrid Boussinesq Saint-Venant model (B-SV), given by the Cauchy problem
\begin{equation}
    \left\{
        \begin{array}{rl}
            \partial_t \eta +   \partial_x u &= 0 \quad \text{ in } Q \\
            \left(1-\chi_B \mu^2 \partial_x^2\right) \partial_t u  +   \partial_x \eta &= 0\quad \text{ in } Q \\ 
            \eta(\cdot, 0) = \eta_0, \quad u(\cdot, 0) &= u_0  \text{ in }\Omega \\ 
        \end{array}\right.
        \label{eq:hybrid-full-adimensional}    
\end{equation}
where $\chi_B$ is the indicator function of $Q_B$, or more compactly as  
\begin{equation}
    \left\{ \begin{split}
        \mathcal{L}_{hyb} \vect W = 0 \quad \text{ in } Q \\ 
        \vect W(\cdot,0) = \vect W_0 \quad \text{ in } \Omega
    \end{split}\right.
    \label{eq:hybrid-compact}
\end{equation}
with $\mathcal{L}_{hyb} = (I - \chi_B A)\partial_t  + B\partial_x$. 
In this case the exchange is not made explicitly through boundary conditions but thanks to the variable coefficient $\chi_B(x)$. The choice of coupling conditions such as \eqref{eq:intro-coupling-conditions} or $\chi$ in \eqref{eq:hybrid-full-adimensional} can drastically change the solution, leading to either different accuracies \cite{GanderMartin2023} and even affecting the well-posedness of the problem (see defective boundary conditions, for example \cite{QUARTERONI2016193}), hence their characterization is crucial for the analysis.

In this work we provide an analysis of system \eqref{eq:hybrid-compact} that could be easily adapted to variations in equations, dimensions, discrete or continuous equations, etc. In \Cref{sec:framework}, we present some facts and tools that will help frame the analysis of \Cref{sec:analysis}. The analysis consists first in deriving the operators $B_\pm$ (as in \eqref{eq:intro-coupling-conditions}) that are implicit in the formulation of the hybrid model \eqref{eq:hybrid-compact}. Then we proceed to derive an analytical solution using harmonic representations (Fourier / Laplace transforms) which we use to prove the Hadamard well-posedness of the problem. This derivation is highly inspired by ideas from domain decomposition \cite{dubach1993contribution}, absorbing boundary conditions \cite{Halpern1982,engquist1977} using similar tools as in \cite{marti2003}. Once the well-posedness is proved we proceed to discuss a way to compute the approximation error (with respect to the full Euler equations) that does not require using the 3D model as a reference, but only a one-way model. We propose a definition for \emph{coupling error} that can be applied to other settings very easily and that explains some of the observations mentioned before. The main advantage being that it isolates the effect that different coupling conditions, such as \eqref{eq:intro-coupling-conditions}, have in the solution. Finally, using the solution that was derived before we show that this coupling error corresponds to the reflections that appear in the solution. Their asymptotic asymptotic size is quantified also. \Cref{sec:conclusions} finishes with the conclusions, and in \cref{sec:appendix} we prove the Hadamard well-posedness of the half-line Boussinesq problem, which extends existing results to a larger class of data.

%% file: figures/intro2/variables.tex
\tikzset{every picture/.style={line width=0.75pt}} %

\begin{tikzpicture}[x=0.75pt,y=0.75pt,yscale=-1,xscale=1]

\draw  [draw opacity=0][fill={rgb, 255:red, 166; green, 198; blue, 230 }  ,fill opacity=1 ] (51.61,35.78) .. controls (52.39,23.12) and (118.6,61.06) .. (136.5,59.43) .. controls (154.4,57.79) and (188.62,22.59) .. (206,21.99) .. controls (223.38,21.4) and (294.09,57.98) .. (313.8,58.31) .. controls (333.5,58.63) and (377.8,33.83) .. (377.8,40.31) .. controls (377.8,46.79) and (377.8,169.66) .. (377.8,178.31) .. controls (377.8,186.95) and (149.96,177.47) .. (143.96,177.47) .. controls (137.96,177.47) and (139.22,179.71) .. (121.96,183.47) .. controls (104.69,187.23) and (52.8,195.01) .. (52.8,189.31) .. controls (52.8,183.6) and (50.83,48.44) .. (51.61,35.78) -- cycle ;
\draw    (240.45,168.33) -- (238.61,49.63) ;
\draw [shift={(238.56,46.63)}, rotate = 89.11] [fill={rgb, 255:red, 0; green, 0; blue, 0 }  ][line width=0.08]  [draw opacity=0] (8.93,-4.29) -- (0,0) -- (8.93,4.29) -- cycle    ;
\draw [shift={(240.5,171.33)}, rotate = 269.11] [fill={rgb, 255:red, 0; green, 0; blue, 0 }  ][line width=0.08]  [draw opacity=0] (8.93,-4.29) -- (0,0) -- (8.93,4.29) -- cycle    ;
\draw    (128.44,116.41) -- (194.16,115.7) ;
\draw [shift={(197.16,115.67)}, rotate = 179.38] [fill={rgb, 255:red, 0; green, 0; blue, 0 }  ][line width=0.08]  [draw opacity=0] (8.93,-4.29) -- (0,0) -- (8.93,4.29) -- cycle    ;
\draw [color={rgb, 255:red, 155; green, 155; blue, 155 }  ,draw opacity=1 ] [dash pattern={on 0.84pt off 2.51pt}]  (33.76,45.87) -- (378.8,47.31) ;
\draw    (206.28,46.59) -- (206.28,24.05) ;
\draw [shift={(206.28,21.05)}, rotate = 90] [fill={rgb, 255:red, 0; green, 0; blue, 0 }  ][line width=0.08]  [draw opacity=0] (8.93,-4.29) -- (0,0) -- (8.93,4.29) -- cycle    ;
\draw [shift={(206.28,46.59)}, rotate = 90] [color={rgb, 255:red, 0; green, 0; blue, 0 }  ][line width=0.75]    (0,5.59) -- (0,-5.59)   ;
\draw  [draw opacity=0][fill={rgb, 255:red, 254; green, 232; blue, 174 }  ,fill opacity=1 ] (52,172.5) -- (379,172.5) -- (379,195.33) -- (52,195.33) -- cycle ;
\draw [color={rgb, 255:red, 0; green, 0; blue, 0 }  ,draw opacity=1 ]   (296.42,22.28) -- (296.62,45.81) ;
\draw [shift={(296.64,47.81)}, rotate = 269.51] [color={rgb, 255:red, 0; green, 0; blue, 0 }  ,draw opacity=1 ][line width=0.75]    (4.37,-1.32) .. controls (2.78,-0.56) and (1.32,-0.12) .. (0,0) .. controls (1.32,0.12) and (2.78,0.56) .. (4.37,1.32)   ;
\draw [shift={(296.41,20.28)}, rotate = 89.51] [color={rgb, 255:red, 0; green, 0; blue, 0 }  ,draw opacity=1 ][line width=0.75]    (4.37,-1.32) .. controls (2.78,-0.56) and (1.32,-0.12) .. (0,0) .. controls (1.32,0.12) and (2.78,0.56) .. (4.37,1.32)   ;
\draw    (52.68,13.88) -- (191.85,13.88) ;
\draw [shift={(193.85,13.88)}, rotate = 180] [color={rgb, 255:red, 0; green, 0; blue, 0 }  ][line width=0.75]    (10.93,-3.29) .. controls (6.95,-1.4) and (3.31,-0.3) .. (0,0) .. controls (3.31,0.3) and (6.95,1.4) .. (10.93,3.29)   ;
\draw [shift={(50.68,13.88)}, rotate = 0] [color={rgb, 255:red, 0; green, 0; blue, 0 }  ][line width=0.75]    (10.93,-3.29) .. controls (6.95,-1.4) and (3.31,-0.3) .. (0,0) .. controls (3.31,0.3) and (6.95,1.4) .. (10.93,3.29)   ;
\draw [color={rgb, 255:red, 155; green, 155; blue, 155 }  ,draw opacity=1 ] (-3.24,200.74) -- (392,200.74)(34.5,17.87) -- (34.5,215.33) (385,195.74) -- (392,200.74) -- (385,205.74) (29.5,24.87) -- (34.5,17.87) -- (39.5,24.87)  ;

\draw (254.38,127.33) node [anchor=north west][inner sep=0.75pt]  [font=\normalsize]  {$h=cte$};
\draw (16,38.33) node [anchor=north west][inner sep=0.75pt]    {$0$};
\draw (131.98,93.33) node [anchor=north west][inner sep=0.75pt]  [font=\normalsize]  {$u( x,t)$};
\draw (210.7,3.71) node [anchor=north west][inner sep=0.75pt]  [font=\normalsize]  {$\eta ( x)$};
\draw (388.8,204.71) node [anchor=north west][inner sep=0.75pt]    {$x$};
\draw (12.8,9.71) node [anchor=north west][inner sep=0.75pt]    {$z$};
\draw (298.41,23.68) node [anchor=north west][inner sep=0.75pt]    {$a$};
\draw (120.53,-3.6) node [anchor=north west][inner sep=0.75pt]    {$L$};

\end{tikzpicture}

%% file: 02-framework.tex
\section{Framework}
\label{sec:framework}
Let $\Omega$ be an open set in $\mathbb{R}$. Let $T>0$ and $Q = \Omega \times ]0,T[$. Let $r,s$ be two non-negative real numbers, most results in this work will be presented in terms of anisotropic spaces of degrees $r\geq 0$ and $s\geq0$ as
$$
    H^{r,s}(Q)  = L^2(0,T; H^r(\Omega)) \cap H^s(0,T; L^2(\Omega))
$$
The analysis will be carried out also using Fourier and Laplace transforms. The direct and inverse Fourier transforms of a function $\phi \in L^2(\mathbb{R})$ are given by
\begin{equation}
    \begin{array}{ll}
         \mathcal{F}\phi(\kappa) = \frac{1}{2\pi} \int_{\mathbb{R}} \phi(x)e^{-j\kappa x} dx, \quad \kappa \in \mathbb{R} \\\\
        \mathcal{F}^{-1}(\mathcal{F}\phi)(x) = \int_{\mathbb{R}} \mathcal{F} \phi(\kappa)e^{j\kappa x} d\kappa, \quad x \in \mathbb{R}
    \end{array}
\end{equation}
where $j$ is the imaginary unit. Also, for $\psi:[0,\infty)\to \mathbb{R}$ such that for some $\alpha>0$, $t \to e^{-\sigma t} \psi(t)$ belongs to $L^2(0,T)$ for any $\sigma>\alpha$,  its Laplace transform is given by
\begin{equation}
    \mathcal{L} \psi(\sigma+j\omega) = \frac{1}{2\pi} \int_0^\infty \psi(x,t) e^{-(\sigma+j\omega)t}dt
\end{equation}
where $s=\sigma+j \omega$ is the Laplace variable with frequency $\omega$ and decay rate $\sigma>\alpha$.  The inverse Laplace transform is then given by
\begin{equation}
    \mathcal{L}^{-1}(\mathcal{L}\psi)(t) = e^{\sigma t}\int_\mathbb{R} \mathcal{L} \psi(\sigma + j \omega) e^{j\omega t}d\omega    
\end{equation}
since $\mathcal L \psi = \mathcal F \phi$ if $\phi = e^{-\sigma t} \psi$ then we can use the usual Fourier-transform theorems such as Parseval's and define the norms 
\begin{equation}
    \begin{array}{l}
        |\phi|_{H^r(\mathbb{R})} = \int (1+|\kappa|^2)^{r/2} \mathcal F \phi(\kappa)\,d\kappa \\
        |\psi|_{H^s(0,\infty)} = \int (1+|\omega|^2)^{s/2} \mathcal L \psi(j\omega)\,d\omega
    \end{array}
\end{equation}
As discussed in \cite{lionsmagenes1961}, a function in $Q \times]0,T[$ can be extended by 0 to $\mathbb{R}^2$ and with these definitions
\begin{equation}
    \begin{array}{ll}
        f\in H^{r,s}(Q) &\Rightarrow |f|_{H^{r,s}(Q)} = \\
        &\int \left[(1+|\kappa|^2)^{r/2} + (1+|\omega|^2)^{s/2} \right] \mathcal L \mathcal F f(\kappa, j\omega)\,d\kappa d\omega < \infty        
    \end{array}
\end{equation}

In the following the "hat notation" $\widehat \phi$ or $\widehat \psi$ will be used to refer to either transform when no confusion is possible, which will be clarified accordingly.

\subsection{Boussinesq equations}
\subsubsection{On the real line}

\label{sec:cauchy-b}
On the real line, the (linear) Boussinesq equations are
\begin{equation}
    \left\{\begin{array}{cc}
        \mathcal{L}_{B} \vect W = 0 & \text{ in } Q \\
        \vect W(0, \cdot ) = \vect W_0    & \text{ in } \Omega
    \end{array}\right.
    \label{eq:lgn-homogeneous-cauchy}
\end{equation}

Denoting by $\widehat \phi$ the Fourier transform of $\phi$, then from \eqref{eq:lgn-homogeneous-cauchy} one gets
\begin{equation}
    \left\{
    \begin{array}{cc}
        Q(j\kappa) \partial_t \widehat{\vect W}(\kappa,t) = P(j\kappa)\widehat{\vect W}(\kappa,t) & (\kappa,t) \in \mathbb{R} \times ]0,T[ \\
        \widehat{\vect W}(\kappa, 0) = \widehat{\vect W}_0(\kappa) &\quad \kappa \in \mathbb{R}
    \end{array}
    \right.
    \label{eq:lgn-homogeneous-cauchy-fourier}
\end{equation}
where $\widehat{\vect W} = (\hat \eta, \hat u)$ and $P(j\kappa) = -j\kappa B$ and $Q(j\kappa) = I + \kappa^2 A $ are the Fourier symbols of the differential operators  $P(\partial_x) = - B\partial_x$ and $Q = (I-A \partial_x^2)$. The solution of \eqref{eq:lgn-homogeneous-cauchy-fourier} can be written as
\begin{equation}
    \widehat{\vect W} = e^{Q^{-1}(j\kappa)P(j\kappa)t} \widehat{\vect W_0}
    \label{eq:lgn-homogeneous-cauchy-exp-solution}
\end{equation}
and taking the inverse Fourier transform one obtains the general formula
\begin{equation}
    \vect{W}(x,t) = \int_\mathbb{R} e^{j\kappa x} e^{Q^{-1}(j\kappa)P(j\kappa) t}\widehat{\vect W}_0(\kappa) d\kappa \quad (x,t) \in Q
    \label{eq:lgn-homogeneous-cauchy-fourier-solution-first}
\end{equation}

Now we can write $Q^{-1}(j\kappa)P(j\kappa) = S J S^{-1}$ with $S$ and $J$ the matrices of eigenvectors and eigenvalues, respectively given by 
\begin{equation}
    S = \left(\begin{array}{cc}
        -\sqrt{1+\mu^2\kappa^2} & \sqrt{1+\mu^2\kappa^2} \\
        1                       & 1
    \end{array}\right)
    \label{eq:b-cauchy-diagonalization}
\end{equation}
\begin{equation}
    J(\kappa) = j\omega(\kappa) \left(
        \begin{array}{cc}
            1 & 0 \\ 
            0 & -1
        \end{array}
    \right), \quad 
    \omega(\kappa) = \frac{\kappa}{\sqrt{1+\mu^2 \kappa^2}}
    \label{eq:dispersion-relation}
\end{equation}
where  $det(S)=-2\sqrt{1+\mu^2\kappa^2}\neq 0$ for any real $\kappa$, so it is well-defined, and $\omega(\kappa)$ is known as the dispersion relation. 

Recalling that $e^{J} = \sum_{k=0}^\infty \frac{J^k}{k!}$, the solution can now be written as 
\begin{equation}
    \widehat {\vect W} =Se^{Jt} S^{-1} \widehat{\vect W}_0 \quad \text{ in } \mathbb{R} \times ]0,T[
\end{equation}
so taking the inverse Fourier transform, we obtain
\begin{equation}
    \vect{W}(x,t) = \int_\mathbb{R} e^{j\kappa x} S(\kappa)e^{J(\kappa)t} S^{-1}(\kappa) \widehat{\vect W}_0(\kappa) d\kappa \quad \text{ in } Q
    \label{eq:lgn-homogeneous-cauchy-fourier-solution}
\end{equation}
and because the Fourier representation is unique \cite{kreiss2004} we have that \cref{eq:lgn-homogeneous-cauchy} has a unique solution given by \eqref{eq:lgn-homogeneous-cauchy-fourier-solution}. 

The Boussinesq equations also admit a global energy conservation law on each harmonic. Let $D(\kappa)$ be 
\begin{equation}
    D(\kappa) = \left(\begin{array}{cc}
        1 & 0                  \\
        0 & 1 + \mu^2 \kappa^2
    \end{array}\right)
\end{equation}
so, defining the complex vector product  $\langle a,b  \rangle = \overline{a}^T b$, one has
\begin{equation}
    \begin{split}
        \frac{d}{dt} \langle \widehat{\vect W}, D \widehat{\vect W}\rangle
        &= \langle \frac{d}{dt}\widehat{\vect W}, D \widehat{\vect W} \rangle + \langle\widehat{\vect W},  D \frac{d}{dt} \widehat{\vect W} \rangle \\
        &= \langle Q^{-1}P \widehat{\vect W}, D \widehat{\vect W} \rangle + \langle\widehat{\vect W},  D Q^{-1}P \widehat{\vect W} \rangle \\
        &=\langle \widehat{\vect W}, \left(Q^{-1}\overline{P}D + DQ^{-1}P \right)\widehat{\vect W}\rangle
    \end{split}
\end{equation}
and
\begin{equation}
    \left(Q^{-1}P\right)^*D = \left(
    \begin{array}{cc}
            0        & j\kappa \\
            j \kappa & 0
        \end{array}
    \right)
    = - DQ^{-1}P
\end{equation}
which means that
\begin{equation}
    \langle \vect{\widehat W}, D \vect{\widehat W}\rangle = \langle \vect{\widehat W}_0, D  \vect{\widehat W}_0\rangle \quad \text{ in } \mathbb{R} \times \mathbb{R}^+
\end{equation}
defining the energy norm of the Boussinesq equation as $|\widehat {\vect W}|_B = \langle \vect{\widehat W}, D \vect{\widehat W}\rangle^{1/2}$ we deduce the following theorem. 
\begin{theorem}
    \label{th:cauchy-b}
    Let $r\geq 0$, if $\vect W_0 =(\eta_0, u_0)\in H^r(\Omega) \times H^{r+1}(\Omega)$ then the Cauchy problem \eqref{eq:lgn-homogeneous-cauchy} has a unique solution $\vect W =(\eta, u) \in L^\infty(0,\infty; H^r(\Omega)\times H^{r+1}(\Omega))$ and

    \begin{equation}
        |(1+\kappa^2)^{r/2}\, \vect{\widehat W}(t, \cdot) |_B = |(1+\kappa^2)^{r/2} \vect{\widehat W}_0 |_B,\quad \forall t \geq 0
        \label{eq:b-energy}
    \end{equation}
\end{theorem}
\begin{remark}
    In the case where $r=n$ is an integer, this can be written as
    \begin{equation}
        \int_\mathbb{R} |\partial_x^n \eta|^2 + |\partial_x^n u|^2  + \mu^2 |\partial_x^{n+1} u|^2 dx = \int_\mathbb{R} |\partial_x^n \eta_0|^2 +|\partial_x^n u_0|^2 + \mu^2 |\partial_x^{n+1}  u_0|^2 dx, \quad \forall t \geq 0
    \end{equation}
\end{remark}

Moreover, since the dispersion relation has compact image:
\begin{equation}
    |\omega(\kappa)| = \frac{1}{\mu} \frac{\mu |\kappa|}{\sqrt{1+\mu^2\kappa^2}}\leq \frac{1}{\mu}, \quad \forall \kappa \in \mathbb{R}
    \label{eq:mu-upperboun}
\end{equation}
the time-spectrum of the solution has compact support, \textit{i.e.}, the solution is smooth in time. The following \Cref{th:b-cauchy-regularity} expresses this fact more precisely.
 
\begin{proposition}
    \label{th:b-cauchy-regularity}
    Let $r \geq 1$, $\vect W_0 \in H^r(Q)\times H^{r+1}(Q)$. If  $\vect W=(\eta,u)$ is the solution of  \eqref{eq:lgn-homogeneous-cauchy} then 
    $$|\vect W|_{H^{r,s}\times H^{r,s}} \leq C|\vect W_0|_{H^{r+1}\times H^{r+1}}$$
    and $\vect W \in \bigcup_{s \geq 0} H^{r,s}(Q) \times H^{r+1,s}(Q)$
\end{proposition}
\begin{proof}
    Taking the time derivative of \eqref{eq:lgn-homogeneous-cauchy-exp-solution} and using the decomposition \eqref{eq:b-cauchy-diagonalization} one obtains 
    \begin{equation}
        \partial_t^s \vect{\widehat W} = (j \omega)^s S R e^{Jt}S^{-1} \vect{\widehat W}_0, \quad R =\left(  
            \begin{array}{cc}
                1 & 0 \\ 
                0 & -1
            \end{array}
        \right)
    \end{equation}
    which means that 
    \begin{equation}
        |\partial_t^s \vect{\widehat W}| \leq |\omega|^s |S| |S^{-1}| |\widehat{\vect W}_0|
    \end{equation}
    letting $\psi(\kappa) = 1+ \mu^2 \kappa^2$ then
    \begin{equation}
        S S^* = \left(\begin{array}{cc}
                2  \psi^2 & 0 \\
                0         & 2
            \end{array}\right), \quad
        S^{-1} (S^{-1})^* = \frac{1}{4}\left(\begin{array}{cc}
                1+\psi^{-2} & 1-\psi^{-2} \\
                1-\psi^{-2} & 1+\psi^{-2}
            \end{array}\right)
    \end{equation}
    so $|S|=\sqrt{2}\psi$ and  $|S^{-1}|=1/\sqrt{2}$, and since $|Re^{Jt}|=1$ we obtain 
    \begin{equation}
        |\partial_t^s \vect{ \widehat W}(\kappa,t)| \leq |\omega(\kappa)|^s (1+\mu^2\kappa^2)^{1/2} |\vect{\widehat W}_0(\kappa)|
        \label{eq:bcauchy-dtdw}
    \end{equation}
    and finally, we can compute
    \begin{equation}
        \begin{split}
            \int_0^T \int_\mathbb{R} |\partial_t^s \vect W(x,t)|^2 dxdt &=
            \int_0^T \int_\mathbb{R} |\partial_t^s  \vect{\widehat W}(\kappa, t)|^2 d\kappa dt  \\ 
            &\leq T \int_\mathbb{R} |\omega(\kappa)|^{2s} | (1+\mu^2 \kappa^2) |\vect{\widehat W}_0(\kappa)|^2 d\kappa, \quad \text{(Eq. \eqref{eq:bcauchy-dtdw})}\\ 
            &\leq T \int_\mathbb{R} \frac{1}{\mu^{2s}} (1+\mu^2 \kappa^2) |\vect{\widehat W}_0(\kappa)|^2d\kappa, \quad \text{(Eq. \eqref{eq:mu-upperboun})} \\ 
            &= \frac{T}{\mu^{2s}}(|\vect W_0|^2 + \mu^2 |\partial_x \vect W_0|^2), \quad \text{(Parseval)}
        \end{split}
    \end{equation}
    Because $r\geq 1$ implies $\vect W_0 \in H^{1}(\Omega)\times H^{1}(\Omega)$ we conclude that $\partial_t^s \vect W \in  (L^2(0,T; L^2(\Omega)))^2$ for every $s \geq 0$.

\end{proof}

\subsubsection{On the half line}
In this section we are interested in the problems on the half lines $\Omega_- = ]-\infty, 0[$ and $\Omega_+ = ]0,\infty[$, $Q_\pm = \Omega_\pm \times ]0,T[$, of functions $\vect W_-: [0,T] \times \overline {\Omega}_-$, $\vect W_+: [0,T] \times \overline {\Omega}_+$ such that
\begin{equation}
    \left\{
    \begin{split}
        \mathcal{L_B}\vect W_\pm= 0    & \quad \text{ in } Q_\pm     \\
        \vect W_\pm(0,\cdot) = 0 & \quad \text{ in } \Omega_\pm \\
        B(u_\pm) = f & \quad \text{ in } \Gamma
    \end{split}\right.
    \label{eq:half-line-problem}
\end{equation}
where $\Gamma=\{(x,t)=(0,t)\}$ and $\vect W_{\Gamma_\pm}=(\eta_{\Gamma_\pm}, u_{\Gamma_\pm})$ is some given data, $B: Q_\pm \to \Gamma$  a trace operator and $f:\Gamma \to \mathbb{R}$ a given source function.

\label{sec:laplace-halfline}
Also, let $\widehat \psi$ denote the Laplace transform of $\psi$ and  
$\vect{\widehat W}_- = (\hat \eta_-, \hat u_-)$. Taking the Laplace transform of \eqref{eq:half-line-problem}
\begin{equation}
    \left\{
    \begin{array}{rl}
        s\hat \eta_- +  \partial_x \hat u_-                          & = 0 \\
        s \hat u_- - s \mu^2\partial_x^2 \hat u_- +  \partial_x \hat \eta_- & = 0
    \end{array}
    \right.
    \label{eq:lgn-laplace}
\end{equation}

The solution of \eqref{eq:half-line-problem}  in $Q_-$ is given by 
\begin{equation}
    \boldsymbol{\widehat W_-} = V e^{\Lambda x} V^{-1}\vect{\widehat W}_{\Gamma } \quad \text{ in } \Omega_{-}    
\end{equation}
with
\begin{equation}
    \lambda(s) = \sqrt{\frac{s^2}{1+\mu^2 s^2}}
    \label{eq:lambda}
\end{equation}
and
\begin{equation}
    \begin{array}{ccc}
        R = V \Lambda V^{-1},
         & V = \left(\begin{array}{cc}
                \frac{ \lambda(s)}{s} & -\frac{\lambda(s)}{s} \\
                1                     & 1
            \end{array}\right),
         & \Lambda(s) = \left(\begin{array}{cc}
                -\lambda(s) & 0          \\
                0           & \lambda(s)
            \end{array}\right)
    \end{array}
    \label{eq:halfline-diagonalization}
\end{equation}
but in order  to satisfy the initial condition, necessarily one must have 
\begin{equation}
    \begin{array}{rlr}
        \vect{\widehat W_-}(x,s) & = \alpha_-(s) \vect v_2(s)  e^{\lambda(s)x } & (s,x) \in \mathbb{C}^+ \times \Omega_-\\
        \vect{\widehat W_+}(x,s) & = \alpha_+(s) \vect v_1(s)  e^{-\lambda(s)x}    & (s,x) \in \mathbb{C}^+ \times \Omega_+
    \end{array}
    \label{eq:half-line-solution}
\end{equation}
for some $\alpha_\pm(s): \mathbb{C}^+ \to \mathbb{C}$. Moreover, the next theorems hold
\begin{theorem}
    Let $u_{\Gamma} \in H^{k}(0,T)$, $k \geq 0$,  there exists a unique solution to problem \eqref{eq:half-line-problem}  given by \eqref{eq:half-line-solution} and  $$
    \vect W \in  \bigcup_{r \in \mathbb{N}_0} H^{r,k+3}(Q_-) \times H^{r,k+2}(Q_-)
    $$ Moreover, there exists $C(T,\mu, r)>0$ such that
    $$
        |\vect  W|_{E^r} \leq C(T,\mu, r) |u_{\Gamma}|_{H^k(0,T)}
    $$
    with $E^r = H^{r,k+1}(Q_-)\times H^{r,k}(Q_-)$
    \label{th:halfplane-boussinesq}
\end{theorem}

Since this problem is related to general transmission problems such as wave generation, and not only to model-coupling, and for the sake of keeping the focus on the coupled models, the reader can find the proof of this theorem  in section \ref{sec:appendix} of the Appendix. %

\subsection{Saint-Venant equations}
We recall that the problem on the real line 
\begin{equation}
    \left\{ \begin{array}{l}
        \mathcal{L}_{SV} \vect W = 0 \quad \text{ in } Q \\ 
        \vect W = \vect W_0 \quad \text{ in } \Omega
    \end{array}\right.   
    \label{eq:sv-cauchy}
\end{equation}
has a unique solution given by 
\begin{equation}
    \begin{split}
        \eta(x,t) = \frac{1}{2}\left(\eta_0(x-t) +\eta_0(x+t)\right) + \frac{1}{2}(u_0(x-t) -u_0(x+t)) \\
        u(x,t) = \frac{1}{2}(\eta_0(x-t) -\eta_0(x+t)) + \frac{1}{2}(u_0(x-t)  + u_0(x+t)) \\
    \end{split}
    \label{eq:sv-cauchy-solution}
\end{equation}
which is equivalet to \eqref{eq:lgn-homogeneous-cauchy-fourier-solution} with $\mu=0$, dispersion relation $\omega(\kappa)= \pm \kappa$ and constant phase speed $c(\kappa)=\pm 1$. 
Since also the solution has the same regularity as the initial data, for later reference, we write this fact as
\begin{proposition}
    \label{th:sv-cauchy}
    If $\eta_0$ and $u_0$ are in $H^r(\Omega)$, with $r\geq 0 $, then problem \eqref{eq:sv-cauchy} has a unique solution $(\eta, u)$ in $(H^{r,r}(\Omega \times ]0,T[))^2 \cap (L^2(0,\infty; H^r(\Omega)))^2$ and the traces $u(x, \cdot)$ are in $H^r(0,T)$ for any $x\in \Omega$.
\end{proposition}

On the real line we have that
\begin{equation}
    \begin{array}{ll}
        \eta_-(x,t) = u_-(x,t) = u_{\Gamma_-}(t+x) H(t+x) & \text{ in } \Omega_-\\ 
        \eta_+(x,t) = u_+(x,t) = u_{\Gamma_+}(t-x) H(t-x) & \text{ in } \Omega_+
        \label{eq:sv-halfline-solution-explicit}
    \end{array}
\end{equation}
with $H=\chi_{ x >0}$ the Heaviside function and
\begin{proposition}
    \label{th:halfline-sv}
    Let $\vect W_{\Gamma_\pm}=(\eta_{\Gamma_\pm}, u_{\Gamma_\pm})$. If  $u_{\Gamma_\pm} \in H^s(0,T)$, $s\geq0$, and $\eta_{\Gamma_\pm} = \pm u_{\Gamma_\pm}$ then there exists a unique solution $\vect W_\pm =(\eta_{\pm}, u_{\pm}) \in H^{s,s}(\Omega_\pm \times ]0,T[)$ given by \eqref{eq:sv-halfline-solution-explicit}, and for every $x \in \Omega_{\pm}$ one has $u_{\pm}(x, \cdot) \in H^s(0,T)$. Moreover, there exists $C(T)>0$ such that
    \begin{equation}
        |\vect W|_{H^{s,s}(\Omega_\pm \times ]0,T[)} \leq C(T) |u_{\Gamma_\pm}|_{H^s(0,T)}
    \end{equation}
\end{proposition}

\pagebreak

\subsection{Useful analysis theorems}

\subsubsection{Trace theorems}
The following theorem from \cite{lionsmagenes1961} characterizes the derivatives of a function on $H^{r,s}(Q)$

\begin{theorem}
\label{th:derivatives-theorem}
    If $u \in H^{r,s}(Q)$, with $r$ and $s$ two strictly positive real numbers, and $j,k$ are two integeres such that $1-\left(\dfrac{j}{r} + \dfrac{k}{s}\right)\geq 0$, then we have
    $$
        \partial_x^j \partial_t^ k u \in H^{\mu,\nu}(Q)
    $$
    with $\dfrac{\mu}{r} = \dfrac{\nu}{s} = 1 - \left(\dfrac{j}{r}+\dfrac{k}{s}\right)$
\end{theorem}

The following theorem is quoted from \cite[Th. 2.1 ]{lionsmagenes1961} and helps characterize a function in $H^{r,s}(Q)$ in terms of its traces and initial condition. 

\begin{theorem}
    Let $r,s > 0$, the mapping
    \begin{equation}
        u \to (f_k, g_j) = \left\{ \begin{array}{c}
            (\partial_t^k u(x,0))_{k<s-1/2, \, s \geq 1/2} \\ \\
            (\partial_x^j u(0,t))_{j<r-1/2, \,r \geq 1/2}
        \end{array}\right.
    \end{equation}
    with $k,j\geq 0$ integers, is continuous from $H^{r,s}(Q)$ to $ H^{p_k}(\Omega) \times H^{\nu_j}(0,\infty)$ where 
    $$
    p_k = \dfrac{r}{s}(s-k - 1/2), \text{ and }\, \dfrac{\nu_j}{s} = \dfrac{r-j-1/2}{r}
    $$
    \label{th:trace}
\end{theorem}
notice that we only quote the continuity of the trace operators and not the surjectivity proved on \cite[Th. 2.3 ]{lionsmagenes1961} that requires compatibility conditions between the traces. This can be avoided thanks to the properties of the Boussinesq equations.

\subsubsection{Approximation theory}
\label{app:approx}
\begin{theorem}
    \label{eq:approx-theory}
    Let $r\geq 0$ integer, if $f \in H^r(\mathbb{R}^+)$ then there exists a family of functions $(f_h)_{h>0}\subset C^\infty$ and positive constants $c,C$ such that  $\hat f_h \in C_0^\infty(\mathbb{R})$ and
    \begin{equation}
        |f^{(n)}-f_h^{(n)}|_{L^2(\mathbb{R})} \leq C h^{r-n} |f^{(r)}|_{L^2(\mathbb{R})}, \quad \hat f=\hat f_h \text{ in } ]\frac{c}{h},\frac{c}{h}[
    \end{equation}
    for $n$ integer between $0$ and $r$
\end{theorem}
\begin{proof}
    Let $\hat{\cdot}$ denote the Fourier transform; $\varphi \in C^\infty(\mathbb{R})$ such that $|\hat \varphi|\leq 1$, $\hat \varphi = 1$ in $(-c,c)$ and denote 
    \begin{equation}
        \varphi_h(x) = \frac{1}{h} \varphi\left(\dfrac{x}{h}\right) \quad, x \in \mathbb{R}
    \end{equation}
    then $\hat \varphi_h = \hat \varphi(h \cdot)$, so $\varphi_h = 1$ in $(-\frac{c}{h}, \frac{c}{h})$, then 
    \begin{equation}
        \begin{split}
            |f^{(n)} - (f * \varphi_h)^{(n)}|^2_{L^2(\mathbb{R})} &= |\omega^n (\hat f - \hat f \hat \varphi_h)|^2_{L^2(\mathbb{R})}
            \\ &= \int_{|\omega|\geq \frac{c}{h} } |\omega|^{2n}|\hat f(\omega) - \hat f(\omega) \hat \varphi_h(\omega) d\omega \\
            &\leq \left(\frac{h|\omega|}{c}\right)^{2r-2n} \int |\omega|^{2n} |\hat f(\omega) - \hat f(\omega) \hat \varphi_h(\omega)|^2 d\omega \\ 
            &\leq 4\left(\frac{ h^{2}}{c^{2}}\right)^{r-n} |f^{(r)}|_{L^2(\mathbb{R})}^2
        \end{split}
    \end{equation}
    wheret the last inequality is obtained from $|\hat \rho_h| < 1$, a triangle inequality and the Fourier's identity for the derivative.
\end{proof}

The following theorem is quoted from \cite{kreiss2004}
\begin{theorem}
    \label{th:unique-extension}
    Let $B_1$ and $B_2$ denoted normed spaces, let $M$ be a dense subspace of $B_1$ and let $B_2$ be complete.  If $S_0: M \to B_2$ is a bounded linear operator, then there is a unique bounded linear operator $S: B_1 \to B_2$ with $Sf = S_0 f$ for all $f \in M$. The operator $S$ is called the extension of $S_0$. 
\end{theorem}

The next theorem shows that smooth functions that integrate to 0 are dense in $L^2(\mathbb{R})$
\begin{theorem}
    \label{th:M-dense-L2}
    Let $M= \{ g \in C_0^\infty(\mathbb{R}): \int g(x) dx =0 \}$ and $f \in H^r(\mathbb{R})$, $r \geq 0$. Then, for every $\varepsilon>0$  there exists a function $g \in M$ such that 
    \begin{equation}
        |f-g|_{H^r(\mathbb{R})} < \varepsilon 
    \end{equation}
\end{theorem}
\begin{proof}
    Because of the density of $C_0^\infty(\mathbb{R})$ in $H^r(\mathbb{R})$ there exists $h\in C_0^\infty$ such that 
    \begin{equation}
        |f - h|_{H^r(\mathbb{R})} < \varepsilon
    \end{equation}
    Let 
    \begin{equation}
        C = \int h(x) dx
    \end{equation}
    and for $L>0$, let $\phi \in C_0^\infty(\mathbb{R})$ a bump function such that $supp(\phi) \subset (-1,1)$ and $\int \phi dx = 1$, and define $\phi_L(x) = \frac{1}{L} \phi(\frac{x}{L})$, then 
    \begin{equation}
        \phi^{(n)}_L(x) = \frac{1}{L^{n+1}} \phi^{(n)}(\frac{x}{L}), \quad \text{ and } \int \phi_L(x)dx = 1
    \end{equation}
    and let $M = \max_n |\phi^{(n)|^2}$, then 
    \begin{equation}
        \begin{split}
            |\phi_L^{(m)}|_{L^2}^2 &= \int_L^{L}|\phi_L^{(n)}(x)|^2 dx  \\ 
            &\leq \frac{1}{L^{2(n+1)}} \int_L^{L}|\phi^{(n)}(\frac{x}{L})|^2 dx  \\ 
            &\leq \frac{M 2L}{L^{2(n+1)}} = \frac{2M}{L^{2n+1}}
        \end{split}
    \end{equation}
    hence 
    \begin{equation}
        |\phi_L|_{H^{n}}\leq \frac{2M(n+1)}{L}
    \end{equation}
    Let $C = \int h dx$ and 
    \begin{equation}
        g(x) = h(x) - C \phi_L(x)
    \end{equation}
    then by construction $g \in M$ and 
    \begin{equation}
        \begin{split}
            |f-g|_{H^r(\mathbb{R})} &\leq |f-h|_{H^r(\mathbb{R})} + C|\phi_L|_{H^{n}} \\ 
            & \leq \varepsilon + \frac{2M(n+1)}{L}
        \end{split}
    \end{equation}
    so choosing $L = 2M(n+1)/\varepsilon$, we conclude that $g\in M$ and 
    \begin{equation}
        |f-g|_{H^r(\mathbb{R})} \leq 2 \varepsilon
    \end{equation}
\end{proof}

\begin{theorem}
    \label{th:M0-dense-L2}
    Let $M= \{ g \in C^\infty(\mathbb{R}): g(0) =0 \}$ and $f \in H^r(\mathbb{R})$, $r \geq 0$. Then, for every $\varepsilon>0$  there exists a function $g \in M$ such that 
    \begin{equation}
        |f-g|_{H^r(\mathbb{R})} < \varepsilon 
    \end{equation}
\end{theorem}
\begin{proof}
    Let $\hat f$ denote the Fourier transform of $f$ and $M_0=\{ g \in C_0^\infty(\mathbb{R}): \int  g(\kappa) d\kappa = 0  \}$, then by \Cref{th:M-dense-L2} there exists $g \in M_0$ such that 
    \begin{equation}
        |\hat f - g|_{H^n} \leq \varepsilon
    \end{equation}
    so denoting $\mathcal{F}^{-1}(g)$ the inverse Fourier transform of $g$ one has $g \in C^\infty$ and 
    \begin{equation}
        g(0) = \int \hat g(\kappa) d\kappa = 0 
    \end{equation}
    so $g \in M$ and from Parseval's theorem
    \begin{equation}
        | f - \mathcal{F}^{-1}(g)|_{H^n} \leq \varepsilon
    \end{equation}
\end{proof}

\begin{theorem}
    \label{th:elliptic}
    Let $n$ be a non-negative integer, $g \in H^n(\mathbb{R}^-)$ and $y_0 \in \mathbb{R}$, then there exists one and only one function $y_1: \mathbb{R}^- \to \mathbb{R}$ such that
    \begin{equation}
        y_1 - \mu^2 y_1'' = g \quad , y(0) = y_0
        \label{eq:elliptic-eq-dirichlet}
    \end{equation}
    Also, there is one and only one function $y_2: \mathbb{R}^- \to \mathbb{R}$ such that
    \begin{equation}
        y_2- \mu^2 y_2'' = g \quad , y'(0) = y_0
    \end{equation}
    also $y_i \in H^{n+2}(\mathbb{R}^-)$ and for integer $k$ from $0$ to $n$
    \begin{equation}
        | y_i^{(k)}|^2_{L^2(\Omega)} + \mu^2 |y_i^{(k+1)}|^2_{L^2(\Omega)} + y^{(k)}(0)y^{(k+1)}(0) \leq \frac{1}{2} |g^{(n)}|^2_{L^2(\Omega)}
    \end{equation}
\end{theorem}
\begin{proof}
The existence and uniqueness are given by the formula
\begin{equation}
    y(x) = \frac{1}{\mu} \int_0^x \sinh((z-x)/\mu) g(z) dz + c e^{x/\mu}
\end{equation}
where $c=y_0$ in the Dirichlet case, and $c=\mu y_0$ in the Neumann case. The smoothness can be obtained by multiplying by $y$ and integrating by parts  to obtain
\begin{equation}
    \int_{-\infty}^0y^2(x) + \mu^2 (y'(x))^2 dx  + y(0)y'(0) - y(-\infty)y'(-\infty) = \int_{-\infty}^0 g(x)y(x) dx
\end{equation}
and using the inequality
\begin{equation}
    ab \leq \frac{a^2}{2\alpha}  + \alpha \frac{b^2}{2}
    \label{eq:ab-inequality}
\end{equation}
with $\alpha=1$ to obtain
\begin{equation}
    \int_{-\infty}^0 y^2 + \mu^2 (y'(x))^2 dx + y(0)y'(0) + y(-\infty)y'(-\infty) \leq \frac{1}{2} \int_{-\infty}^0 g(x)^2 dx
    \label{eq:h1-estimate}
\end{equation}
which implies that $y \in H^1(\mathbb{R}^-)$, and $y(-\infty)=y'(-\infty)=0$. For the regularity, because of the ODE, we have $y''=(g-y)/\mu \in L^2(\mathbb{R}^-)$ so $y \in H^{2}(\mathbb{R}^+)$. Now, for integer $k$ from $1$ to $n$ by linearity  $\tilde y = y^{(k)}$ and $\tilde g = g^{(k)}$ also satisfy \eqref{eq:elliptic-eq-dirichlet}, from where one obtains the energy estimates and $y \in H^{n+2}(\mathbb{R})$.
\end{proof}

%% file: 03-hybrid-analysis.tex
\section{Analysis of  coupled models}
\label{sec:analysis}

To analyze \cref{eq:hybrid-compact} we will derive an analytical solution by writing
\begin{equation}
    \begin{split}
        \vect W = \vect W^\star+ \vect W' & \quad \text{ in } Q
        \label{eq:exact-solution-abstract}
    \end{split}
\end{equation}
with $\vect W^\star$ a convenient known reference, and then solving for $\vect W'$ instead. In the following we propose that choosing $\vect W^\star$ as the one-way model is a simple-yet-powerful choice for several reasons. The first of them is that $\vect W'$ can be understood as reflected waves that propagate both-ways from the interface. In fact, we will solve for $\vect W'$ using a reflection method known for initial-boundary value problems of the 1D wave equation \cite{evans1988}. An important distinction is that in our case the amplitude of the reflection, the reflection coefficient, depends on the frequency, \textit{i.e.}, it is a filter. Such a reflection coefficient has also been derived in the study of absorbing boundary conditions \cite{engquist1979rad,engquist1977,GIVOLI2004319, Halpern1982}. To the authors' knowledge, this is the first time it is used to study the coupling of wave equations with different dispersion relations.

A second reason why the choice of the one-way model for the reference solution $\vect W^\star$ is convenient, is that it facilitates the estimation of the error 
\begin{equation}
    \vect W - \vect W_{3D}
\end{equation}
where $\vect W_{3D}$ is the solution of the free-surface Euler equations, or any other expensive and more accurate model not necessarily involved in the coupling. In fact, we will discuss that the size of $\vect W'$ readily quantifies the contribution that the coupling conditions add in isolation and without relying on $\vect W_{3D}$. This can help distinguish which coupling conditions are better, such as boundary conditions $B_\pm$ in \eqref{eq:intro-coupling-conditions}, or the heterogeneous coefficient $\chi$ in \eqref{eq:hybrid-full-adimensional}, or another. Al

\subsection{The reference solution}

Assuming that the initial data $\vect W_0$ is supported on $\Omega_-$ and that either $\mathcal{L}_- = \mathcal{L}_B$ and $\mathcal{L}_+=\mathcal{L}_{SV}$, the BSV case, or $\mathcal{L}_- = \mathcal{L}_{SV}$ and $\mathcal{L}_+=\mathcal{L}_{B}$, the SVB case, the \emph{one-way} model can be written as 

\begin{equation}
    \vect W^\star = \left\{ 
    \begin{array}{cc}
        \vect W_{cauchy,-}^\star & \text{ in } Q_- \\ 
        \vect W_{half,+}^\star & \text{ in } Q_+
    \end{array}
    \right.
    \label{eq:one-way-model-right-moving}
\end{equation} 
where
\begin{equation}
    \left\{ 
        \begin{array}{cc}
            \mathcal{L}_{-} \vect W_{cauchy,-}^\star =0, & \quad \text{ in } Q_- \\ 
            \vect W_{cauchy,-}^\star (\cdot, 0) = \vect W_0 & \quad \text{ in } \Omega_-
        \end{array}
    \right.
\end{equation}
and
\begin{equation} 
    \left\{ 
        \begin{array}{cc}
            \mathcal{L}_{+} \vect W_{half,+}^\star =0, & \quad \text{ in } Q_+ \\ 
            \vect W_{half,+}^\star(\cdot, 0) = 0 & \quad \text{ in } \Omega_+ \\ 
            u_{half,+}^\star(0,\cdot) = u_{cauchy,-}^\star(0,\cdot) & \quad \text{ on } ]0,T[
        \end{array}
        \label{eq:one-way-model-half}
    \right., 
\end{equation}
Later on, the assumption that $\vect W_0$ is supported in $\Omega_-$ will be also relaxed by using linear superposition of two one-way models. 
Let us examine a few of its properties now.  Substituting solutions of the form 
\begin{equation}
    u_{cauchy,-}^\star(x,t) = e^{j(\kappa_-x  - \omega_-(\kappa_-)t)}, \quad u_{half,+}^\star(x,t) = e^{j(\kappa_+ x - \omega_+(\kappa_+)t)}
\end{equation}
with positive $\kappa_\pm$, and $\omega_\pm(\kappa_\pm)$ the dispersion relations associated to $\mathcal{L}_\pm$, the one-way transmission condition at $x=0$ of \eqref{eq:one-way-model-half} implies that 
\begin{equation}
    \omega_-(\kappa_-) = \omega_+(\kappa_+)
    \label{eq:oneway-compatibility}
\end{equation}
This is a compatibility condition that ensures that the solution of the one-way model is continuous across the interface.

In the case where $Q_-=Q_B$ and $Q_+ = Q_{SV}$ \cref{eq:oneway-compatibility}  can be written as 
\begin{equation}
    \kappa_+ = \frac{\kappa_-}{\sqrt{1+\mu^2\kappa_-^2}} < \kappa_-
\end{equation}
which reveals that when a wave crosses the interface and arrives into $\Omega_+$ its wave length increases. This is an artificial refraction induced by the change in the dispersion relation of the equations. On the opposite direction, when $Q_-=Q_{SV}$ and $Q_+=Q_B$, there is also an artificial refraction and, for the same reason as before, the wave lengths decrease according to
\begin{equation}
    \kappa_+ = \dfrac{\kappa_-}{\sqrt{1-\mu^2\kappa_-^2}} > \kappa_-,\quad \text{ if } \mu \kappa_-<1
\end{equation}
but in the case $\mu \kappa_- >1$ the mapping
\begin{equation}
    \begin{array}{rrll}
        \mathbb{R}\to& 
        ]-\frac{1}{\mu},&\frac{1}{\mu}[&\to 
        \mathbb{R} \\
        \kappa_- \mapsto& 
        \omega_-(\kappa_-) =& 
        \omega_+(\kappa_+) & 
        \mapsto \kappa_+
    \end{array}
\end{equation}
does not work anymore, unless also complex wave numbers are allowed. This issue is associated to a smoothing effect of the half-line (wave generation) problem of the Boussinesq equation studied in \cite[Ch. 6]{galaz2024coupling}.

\subsection{Characterization of the coupling conditions}

The strategy for the analysis of the variable-coefficient Cauchy problem \eqref{eq:hybrid-compact} is to reformulate it as two coupled homogeneous problems in the half-line:
\begin{equation}
    \left\{ \begin{array}{ll}
        \mathcal{L}_- \vect W_- = 0, &\quad \text{ in } Q_- \\ 
        \vect W_-(\cdot,0) = \vect W_0, &\quad \text{ in } \Omega_-
    \end{array}\right. \\ 
    \left\{ \begin{array}{ll}
        \mathcal{L}_+ \vect W_+ = 0, &\quad \text{ in } Q_+ \\ 
        \vect W_+(\cdot,0) = \vect W_0, &\quad \text{ in } \Omega_+
    \end{array}\right.
    \label{eq:hybrid-halflines}
\end{equation}
by finding coupling conditions such that the solution is preserved, \text{i.e.},  
\begin{equation}
    \begin{split}
        \vect W = \vect W_-, & \quad \text{ in } Q_- \\ 
        \vect W = \vect W_+, & \quad \text{ in } Q_+
    \end{split}
\end{equation}
with $\vect W$ solution to \eqref{eq:hybrid-compact} and $(\vect W_-, \vect W_+)$ solution of \eqref{eq:hybrid-halflines}. 
The next theorem reveals the coupling conditions that reformulate problem \eqref{eq:hybrid-compact} as an equivalent coupling of the two half-line problems \eqref{eq:hybrid-halflines}.

\begin{theorem}[Coupling conditions of the hybrid model]

    \label{th:hybrid-coupling-conditions}
    Suppose $\vect W=(\eta, u)$ is solution of \eqref{eq:hybrid-compact}, $\vect W_-=(\eta_-, u_-)$, $\vect W_+=(\eta_+, u_+)$ solve problems \eqref{eq:hybrid-halflines}, and the Laplace transforms of $u_-(,\cdot)$, $u_+(0,\cdot)$, $\partial_x u_-(0,\cdot)$, $\partial_x u_+(0,\cdot)$ exist. Then $\vect W_- = \vect W_{| \overline Q_-}$ and $\vect W_+ = \vect W_{| \overline{Q}_+}$ if and only if
\begin{equation}
    \begin{split}
        u_-(0,\cdot) = u_+(0,\cdot) \\ 
        \partial_x u_-(0,\cdot) = \partial_x u_+(0,\cdot)
    \end{split}
    \quad \text{ in }]0,T[
    \label{eq:coupling-conditions}
\end{equation}
\end{theorem}

\begin{proof}
     Let 
     $$\vect V_- = (\zeta_-,v_-) = \vect W_- - \vect W \quad \text{ in } \overline Q_-$$ 
     $$\vect V_+ = (\zeta_+,v_+) = \vect W_+ - \vect W \quad \text{ in} \overline Q_+$$ 
     then $\vect V_-$ and $\vect V_+$ satisfy \eqref{eq:hybrid-halflines} with 
     \begin{equation}
        \vect V_\pm(\cdot,0) = 0, \quad \text{ in }  \Omega_\pm
     \end{equation}
     
     ($\Rightarrow$): Assuming that $\vect W_- = \vect W_{| \overline Q_-}$ and $\vect W_+ = \vect W_{| \overline{Q}_+}$, then $v_-$ and $v_+$ are constant and equal to zero in $\overline{Q}_\pm$ as also are $\partial_x v_\pm$. Since both $\overline Q_\pm$ contain $\{x=0,t\in[0,T]\}$, we have that  $v_-(0,t)=v_+(0,t)=0$ and $\partial_x v_-(0,t)=\partial_x v_+(0,t)=0$. From the assumption that the Laplace transforms of $u_-(,\cdot), u_+(0,\cdot), \partial_x u_-(0,\cdot), \partial_x u_+(0,\cdot)$ exist, we deduce $\partial_x v_\pm(0,\cdot) = \partial_x u_\pm(0,\cdot) - \partial_x u(0,\cdot)=0$, hence \cref{eq:coupling-conditions} holds. 
   
    ($\Leftarrow)$:  On the other side, let us denote by $\hat{()}$ the Laplace transform. Let us also assume that the equations \eqref{eq:coupling-conditions} hold and let $\lambda_\pm(s)$ be such that 
    \begin{equation}
        \begin{split}
            v_-(x,s) &= \hat u_-(0,s) e^{\lambda_-(s)x},\quad (x,s) \in \Omega_- \times \mathbb{C}^+ \\
            v_+(x,s) &= \hat u_+(0,s)  e^{-\lambda_+(s)x},\quad (x,s) \in \Omega_+ \times \mathbb{C}^+ 
            \label{eq:half-hybrid-proof-solution}
        \end{split}
    \end{equation}
    as in \Cref{sec:laplace-halfline}. Because at $t=0$ $v_-$ and $v_+$ are null, from \eqref{eq:half-hybrid-proof-solution} we know that
    \begin{equation*}
        \begin{split}
            \partial_x \hat v_-(0,s) &= \lambda_-(s) \hat v_-(0,s)\\
            \partial_x \hat v_+(0,s) &= -\lambda_+(s) \hat v_+(0,s)                
        \end{split}
        \label{eq:abstract-transmission-conditions}
    \end{equation*}
    Substracting both equations in \eqref{eq:abstract-transmission-conditions} 
    \begin{equation*}
        \partial_x \hat v_-(0,s) - \partial_x \hat v_+(0,s) = (\lambda_-(s) \hat v_-(0,s) + \lambda_+(s) \hat v_+(0,s))
    \end{equation*}
    and, by definition of $v_\pm$, if the first equation of \eqref{eq:coupling-conditions} holds then
    \begin{equation}
        \partial_x \hat u_-(0,s) - \partial_x \hat u_+(0,s) = (\lambda_-(s)+\lambda_+(s)) 2 \hat v_i(0,s)
        \label{eq:coupling-conditions-equivalence-proof}
    \end{equation}
    is true for both $\hat v_i = \hat v_-$ and $\hat v_i=\hat v_+$. Then, if the second equation of \eqref{eq:coupling-conditions} also holds, the left hand side of \eqref{eq:coupling-conditions-equivalence-proof} is zero. Since $Re(\lambda_-)>0$ and $Re(\lambda_+)>0$, necessarily $v_+$ and $v_-$ are null. Then using the definition of $v_\pm$ we obtain the conclusion.
\end{proof}

At the same time, the coupling conditions \eqref{eq:coupling-conditions} can be reformulated in multiple equivalent ways.  In particular, if instead one defines transmission conditions using linear operators $B_-$, $B_+$ as 
\begin{equation}
    \begin{split}
        B_-(u_-(0,\cdot)) = B_-(u_+(0,\cdot))  \\
        B_+(u_-(0,\cdot)) = B_+(u_+(0,\cdot))  \\
    \end{split}
    \quad \text{ in } ]0,T[
    \label{eq:abstract-bis}
\end{equation} 
instead of \eqref{eq:coupling-conditions}, to preserve the solution of the hybrid model it is enough to choose $B_-$ and $B_+$ such that their Laplace transforms satisfy 
\begin{equation}
    \begin{split}
        \widehat{B u_-(0,\cdot)}   &= 
        \partial_x \widehat{u}_-
        + m_- \hat u_- + c_- \\ 
        \widehat{B u_+(0,\cdot)}   &= 
        \partial_x \widehat{u}_+
        + m_+ \hat u_+ + c_+         
    \end{split}
    \quad \text{ in } \mathbb{C}^+
    \label{eq:coupling-conditions-equivalence}
\end{equation}
with $m_\pm(s), c_\pm(s): \mathbb{C}^+ \to \mathbb{C}$ such that $m_- - m_+ \neq 0$, which is the condition for equations \eqref{eq:abstract-bis} to be linearly independent in the frequency domain of the Laplace transform. This "trick" was noted by \cite{lions1990schwarz} in the context of domain decomposition methods and nowadays is commonly used to define optimized Schwarz methods \cite{gander2006optimized}.

\subsection{Analytical solution}

Let us  denote by $(\lambda_-(s), \vect v_1^-, \vect v_2^-)$ and $(\lambda_+(s)$, $\vect v_1^+$, $\vect v_2^+)$ the positive eigenvalue and pair of eigenvectors of the diagonalization \eqref{eq:halfline-diagonalization}, for  $\mathcal{L}_-$ and $\mathcal{L}_+$ respectively, which is valid for both $\mu=0$, the SV case, and $\mu>0$, the B case. 

To derive the solution for the perturbation $\vect W'$ we will use the transmission conditions 
\begin{equation}
    \left\{ \begin{array}{rll}
        (\partial_x  + \lambda_+) \hat u_-(0,\cdot) 
        &= (\partial_x  + \lambda_+ )\hat u_+(0,\cdot) 
          &
        \quad \text{ in } \mathbb{C}^+ \\ 
        u_+(0,\cdot) 
        &= u_-(0,\cdot)  
        &\quad \text{ in } ]0,T[
    \end{array}\right.
    \label{eq:solution-transmission-conditions}
\end{equation}
which satisfy \eqref{eq:coupling-conditions-equivalence}, so they are equivalent to \eqref{eq:coupling-conditions} and by Theorem \eqref{th:hybrid-coupling-conditions} the solution of \eqref{eq:hybrid-halflines} will be the same as the solution of the hybrid model. Moreover, because $supp(\vect W_0)\subset \Omega_-$, like in \cref{eq:half-line-solution}, $\vect W_+$ satisfies 
\begin{equation}
    \vect{\widehat W}_+(x,s) = \hat u(0,s) \vect v_1^+ e^{-\lambda_+(s)x},\quad (x,s) \in \Omega_+\times \mathbb{C}^+
\end{equation}
and the right hand side of the first line of equation \eqref{eq:solution-transmission-conditions} is 
\begin{equation}
    (\partial_x  + \lambda_+ )\hat u_+(0,\cdot)  = 0 
\end{equation}
which is the exact absorbing boundary condition of the operator $\mathcal{L}_+$. In the BSV the operator is local and corresponds to the transport operator
\begin{equation}
    (\partial_t + \partial_x)u_-(0,\cdot) = 0 
\end{equation}
while in the SVB case it is a pseudo-differential operator and was described in \cite{kazakova2020discrete}.

Because of \eqref{eq:solution-transmission-conditions}, denoting $\vect W_-'=(\eta_-',u_-') = \vect W_- - \vect W^{\star}_{cauchy,-}$ :
\begin{equation}
    (\partial_x + \lambda_+) \hat u_-'(0,\cdot)= -(\partial_x + \lambda_+) \hat{ u}^{\star}_{cauchy,-}(0,\cdot) \quad \text{, in } \mathbb{C}^+
    \label{eq:transmission-conditions-perturbation}
\end{equation}
and because $\hat{u}^{\star}_{cauchy,-}$ also satisfies the half-line problem \eqref{eq:half-line-solution} in $Q_+$ with trace $u^\star_{cauchy,-}(\cdot,)$, we have that
\begin{equation}
    \begin{array}{rll}
        \partial_x \widehat{u}^{\star}_{cauchy,-} 
            &= -\lambda_- \widehat{u}^{\star}_{cauchy,-} 
            & \text{ in } \overline{\Omega_+} \times \mathbb{C}^+ \\ 
        \partial_x \widehat u'_- 
        &= \lambda_- \widehat u'_- 
        & \text{ in } \overline{\Omega_-} \times \mathbb{C}^+ 
    \end{array}
\end{equation}
substituting into \eqref{eq:transmission-conditions-perturbation} we obtain that
\begin{equation}
    \hat u_-'(0,s) = r(s) \widehat{u}^{\star}_{cauchy,-}(0,s), \quad s \in \mathbb{C}^+ \quad \text{ with }\quad  r(s) = \frac{\lambda_-(s)-\lambda_+(s)}{\lambda_-(s)+\lambda_+(s)}
    \label{eq:reflection-coefficient-appearance}
\end{equation}
from \eqref{eq:reflection-coefficient-appearance} we can obtain $\vect W_-'$ and $\vect W_+'$ from the propagation of the signal 
\begin{equation}
    u_\Gamma = \mathcal{L}^{-1}(r \hat u^\star_{cauchy,-}(0,\cdot) ) \quad \text{ in }]0,T[
\end{equation}
on each half-line problem \eqref{eq:half-line-problem}, however, a more expressive formula can be obtained. Using \eqref{eq:half-line-solution} for $\vect W_-'$ and $\vect{W}^{\star}_{cauchy,-}$, we know that
\begin{equation}
    \begin{array}{rll}
        \vect{\widehat W}_-'(x,s) 
        &= r(s) \widehat{ u}^{\star}_{cauchy,-}(0,s) \vect v_2^- e^{\lambda_-(s)x} \quad 
        & \text{ in } \Omega_- \times \mathbb{C}^+ \\ 
        \vect{\widehat{W}}^{\star}_{cauchy,-}(x,s) 
        &=  \widehat{ u}^{\star}_{cauchy,-}(0,s) \vect v_1^- e^{-\lambda_-(s)x} \quad 
        & \text{ in } \Omega_+ \times \mathbb{C}^+
    \end{array}
    \label{eq:reflection-laplace-interm-formulas}
\end{equation}
and since 
\begin{equation}
    \vect v_2^- = G \vect v_1^- \quad \text{ with } G = \left(\begin{array}{cc}
        -1 & 0 \\
        0 & 1
    \end{array}\right)
    \label{eq:g-definition}
\end{equation}
we obtain
\begin{equation}
    \vect{\widehat W}_-'(x,s) = r(s) G \vect{\widehat{ W}}^{\star}_{cauchy,-}(-x,s) \quad \text{ in } \overline{\Omega}_- \times \mathbb{C}^+
    \label{eq:reflection-omega1-laplace}
\end{equation}
This expression shows that the perturbation of the one-way coupling solution can be described by reflecting $\vect W_{cauchy,-}^\star$ through the $x$ axis and by filtering the solution in time by $r(s)$. This is why we refer to $\vect W_\pm'$ as reflections, even though $\vect W_+'$ is in the opposite direction. The formulation \eqref{eq:reflection-omega1-laplace} of the reflections may seem (anti) symmetric because $r(s)$ in \eqref{eq:reflection-coefficient-appearance} only changes sign when $\mathcal{L}_-$ is $\mathcal{L}_B$ or $\mathcal{L}_{SV}$. By using a reflection technique similar to the one sometimes used for the 1D wave equation \cite{evans1988}, the following sections show a more explicit representation that reveals the differences between these two cases.

\begin{remark}
    The reflection coefficient could be derived in the BSV case by noticing that \eqref{eq:solution-transmission-conditions} can be written as 
    \begin{equation}
        (\partial_t + \partial_x)u_-(0,\cdot) = 0 
    \end{equation}
    so, following a classical procedure on absorbing boundary conditions \cite{engquist1977,engquist1979rad,Halpern1982} one can replace $u_-(x,t) = e^{j (\kappa x - \omega(\kappa)t)}+r(\kappa)e^{j (-\kappa x - \omega(\kappa)t)}$, with $\omega(\kappa)$ the dispersion relation of the Boussinesq model, evaluate at $x=0$ and obtain
    \begin{equation}
        r(\kappa) = \frac{\sqrt{1+\mu^2\kappa^2}-1}{\sqrt{1+\mu^2\kappa^2}+1}
    \end{equation}
    which coincides with \eqref{th:b-reflection}. However, in the SVB case, the transparent boundary condition \eqref{eq:solution-transmission-conditions} is a nonlocal pseudo-differential operator \cite{kazakova2020discrete}, and the application of this simple method is not clear. With the current derivation, not only this case can be handled, but also a precise analysis of the model can be derived for arbitrary initial data, which is not possible with the more classical approach.
\end{remark}

\subsubsection{Boussinesq to Saint-Venant case}
We start with the case $\mathcal{L}_- = \mathcal{L}_{B}$ and $\mathcal{L}_+ = \mathcal{L}_{SV}$ (BSV) using the dispersion relation in \eqref{eq:dispersion-relation} that will denoted by $\omega_-(\kappa)$. In the next \Cref{th:b-reflection} the Taylor expansion of $r(s)$  with convergence region $|\mu s|<1$ 
\begin{equation}
    r_N(s) = \sum_{n=1}^N a_n (\mu s)^{2n}
    \label{eq:truncated-coefficient}
\end{equation}
will simplify the inversion of \eqref{eq:reflection-omega1-laplace} and reformulate $\vect W_-'$ as the solution of a Cauchy problem.
\begin{theorem}
    \label{th:b-reflection}
    Let $\mathcal{L}_-=\mathcal{L}_B$. If $\vect W_0$ is supported on $\Omega_-$ then $\widehat{\vect W}_-'$ is the restriction to $Q_-$ of the solution \cref{eq:lgn-homogeneous-cauchy-fourier-solution} of the Cauchy problem  
    \begin{equation}
        \left\{\begin{array}{cc}
            \mathcal{L}_{B} \vect W'_- = 0 & \text{ in } Q \\
            \vect W'_-(0, \cdot) = \mathcal{F}^{-1}\left(r(j\omega_-)\right) * G R \vect W_0   & \text{ in }  \Omega
        \end{array}\right.
    \end{equation}  
    where $R$ is such that $Rf(x) = f(-x)$ and $\mathcal{F}^{-1}$ is the inverse Fourier transform.
\end{theorem}

\begin{proof}
    Let $ \vect{\widehat W}^N$ be such that 
    \begin{equation}
        \widehat{\vect W}^N(x,s) = r_N(s) G \vect{\widehat{ W}}^\star_-(-x,s) \quad x \in \Omega_-
    \end{equation}
    taking the inverse Laplace transform this means that 
    \begin{equation}
        \vect W_-^N(x,t) = \mathcal{L}^{-1}(r_N(s) )* G \vect{ W}^\star_-(-x,\cdot)
    \end{equation}
    since it is a finite sum, we can write 
    \begin{equation}
        \mathcal{L}^{-1}(r_N(s)) * f = \sum_{n=1}^N a_n \mu^{2n} \partial_t^{2n}f 
    \end{equation}
    also, from  \eqref{eq:lgn-homogeneous-cauchy-fourier-solution} 
    \begin{equation}
        \vect{W}(-x,t) = \int_\mathbb{R} e^{j\kappa x} S^{-1}e^{Jt} S \widehat{\vect W}_0(-\kappa) d\kappa
    \end{equation}    
    which means that
    \begin{equation}
        \begin{split}
            \vect{W}_-^N(x,t) &= \mathcal{L}^{-1}(r_N(s) )* G \vect{ W}^\star_-(-x,\cdot) \\ 
            &= \sum_{n=1}^N a_n \mu^{2n}\partial_t^{2n} G \vect{ W}^\star_-(-x,\cdot)\\
            &= \int_\mathbb{R}\sum_{n=1}^N a_n \mu^{2n}  e^{j\kappa x} S^{-1}\partial_t^{2n}\left(e^{Jt}\right) S G \widehat{\vect W}_0(-\kappa) d\kappa
        \end{split}
    \end{equation}
    and since 
    \begin{equation}
        \partial_t^{2n} e^{Jt} = (j \omega_-(\kappa))^{2n} e^{Jt}
    \end{equation}
    after rearranging one can substitute the definition of $r_N$ as
    \begin{equation}
        \begin{split}
            \vect{W}_-^N(x,t) &= \int_\mathbb{R}  e^{j\kappa x} S^{-1}e^{Jt} S \left(\sum_{n=1}^N a_n (\mu j\omega)^{2n} G \widehat{\vect W}_0(-\kappa)\right) d\kappa \\ 
            &= \int_\mathbb{R}  e^{j\kappa x} S^{-1}e^{Jt} S \left(r_N(j\omega) G \widehat{\vect W}_0(-\kappa)\right) d\kappa             
        \end{split}
    \end{equation}   
    since $\mu|\omega(\kappa)|<1$ we can take the limit $N \to \infty$ to obtain
    \begin{equation}
        \vect{W}_-'(x,t) = \int_\mathbb{R}  e^{j\kappa x} S^{-1}e^{Jt} S \left(r(j\omega) G \widehat{\vect W}_0(-\kappa)\right) d\kappa
    \end{equation}
    and comparing with \eqref{eq:lgn-homogeneous-cauchy-fourier-solution} we conclude that $\vect W_-'$ is the solution of \eqref{eq:lgn-homogeneous-cauchy} with initial condition $\mathcal{F}^{-1}\left(r(j\omega)  \right) * G R \vect W_0$.
\end{proof}

\subsubsection{Saint-Venant to Boussinesq case}

The key argument for the convergence of the integrals involved in the proof of \Cref{th:b-reflection} is that the dispersion relation satisfies $\mu |\omega|<1$, so the Taylor series converges. Under this assumption one can repeat the same argument \eqref{th:b-reflection} to deduce the following Proposition.

\begin{proposition}
    If $\mathcal{L}_-=\mathcal{L}_{SV}$, $\vect W_0$ is supported on $\Omega_-$ and $\mathcal{F}(\vect{W}_0)$ is supported on $]-\dfrac{1}{\mu},\dfrac{1}{\mu}[$, then $\vect{ W}'_-$ is the restriction to $\overline{Q}_-$ of the solution \cref{eq:lgn-homogeneous-cauchy-fourier-solution} of the Cauchy problem  
    \begin{equation}
        \left\{\begin{array}{cc}
            \mathcal{L}_{SV} \vect W'_- = 0 & \text{ in } Q \\
            \vect W'_-(0, \cdot ) = \mathcal{F}^{-1}\left(r(j\omega_-)\right) * G R \vect W_0    & \text{ in }  \Omega
        \end{array}\right.
    \end{equation}    
\end{proposition}
The assumption on the spectrum is too strong, however, using formula \eqref{eq:sv-cauchy-solution} and the assumption on the support of the initial data one obtains the following:

\begin{theorem}
    \label{th:sv-reflection}
    If $\mathcal{L}_-=\mathcal{L}_{SV}$, $\vect W_0$ is supported on $\Omega_-$ then $\vect{ W_-}'=(\eta_-', u_-')$ is the restriction to $\overline{Q}_-$ of the solution of \cref{eq:sv-cauchy} of the Cauchy problem  
    \begin{equation}
        \left\{\begin{array}{cc}
            \mathcal{L}_{SV} \vect W'_- = 0 & \text{ in } Q \\
            \vect W'_-(0,\cdot) = H \mathcal{L}^{-1}(r) * G R \vect W_0 & \text{ in } \Omega 
        \end{array}\right.
    \end{equation}    
\end{theorem}

\begin{proof}
    From \eqref{eq:sv-cauchy-solution}, because $\vect W_0 =0$ in $Q_+$
\begin{equation}
    {u}^\star_{cauchy,-}(0,t) = \frac{1}{2}(\eta_0+u_0)(-t) \quad t\geq 0
\end{equation}
substituting on the first equation of \eqref{eq:reflection-laplace-interm-formulas} we obtain 
\begin{equation}
    \vect{\widehat W}_-'(x,s)=\frac{r(s)}{2} (\widehat{R\eta_0}(s) + \widehat{Ru_0}(s)) \left( \begin{array}{c}
        -1 \\
        1
    \end{array}\right) e^{sx}
\end{equation}
so taking the inverse Laplace transform  
\begin{equation}
    \begin{split}
        \eta_-'(x,t) = - \mathcal{L}^{-1}\left(\frac{r}{2} (\widehat{R\eta_0} + \widehat{Ru_0}) \right) ( t+x) H(t+x)  \\ 
        u_-'(x,t) = \mathcal{L}^{-1}\left(\frac{r}{2} (\widehat{R\eta_0} + \widehat{Ru_0}) \right) ( t+x) H(t+x) 
    \end{split}
    ,\quad \text{ in } \Omega_-
\end{equation}
if we now denote $ \eta'_0(x) =  - \mathcal{L}^{-1}(r \widehat{R\eta_0} )$ and $' u_0(x) =  \mathcal{L}^{-1}(r \widehat{R u_0})$ then 
\begin{equation}
    \begin{split}
        \eta_-'(x,t) &= \frac{1}{2}( \eta'_0(x+t) -  u'_0(x+t))\\ 
        u_-'(x,t) &= -\frac{1}{2}( \eta'_0(x+t) -  u'_0(x+t))
    \end{split}
\end{equation}
which contrasted with \eqref{eq:sv-cauchy-solution} proves the conclusion.
\end{proof}
\begin{remark}
    \label{remark:hybrid-solution}
    We can now fully reconstruct the solution of the hybrid model
    \begin{equation}
        \vect W = \left\{
            \begin{array}{rl}
                \vect W_- & \text{ in } Q_- \\
                \vect W_+ & \text{ in } Q_+
            \end{array}
         \right.
    \end{equation}
    In the BSV case we have
    \begin{equation}
        \left\{
        \begin{split}
            \vect W_-(x,t) &= \int_\mathbb{R} e^{j\kappa x} S(\kappa)e^{J(\kappa)t} S^{-1}(\kappa)  \vect W_0^- \,d\kappa \\
            \vect W_+(x,t) &= u_-(0,\cdot) H(t-x) \left(
                \begin{array}{c}
                    1 \\
                    1
            \end{array}\right)
        \end{split} \right.
    \end{equation}
    with 
    \begin{equation}
        \mathcal{F}({\vect W}_0^- ) = (I + rGR)\mathcal{F}(\vect W_0)(\kappa)\, , \kappa \in \mathbb{R}
        \label{eq:hybrid-bsv-solution-ic}
    \end{equation}
    while in th SVB case
    \begin{equation}
        \left\{
        \begin{split}
            \eta_-(x,t) &= \frac{1}{2}\left( \eta_0^-(x-t) + \eta_0^-(x+t)\right) + \frac{1}{2}( u_0^-(x-t) - u_0^-(x+t)) \\
            & \quad\quad\quad(x,t) \in Q_- \\
            u_-(x,t) &= \frac{1}{2}( \eta_0^-(x-t) - \eta_0^-(x+t)) + \frac{1}{2}( u_0^-(x-t)  +  u_0^-(x+t)) \\
            & \quad\quad\quad(x,t) \in Q_- \\
            \mathcal{L}(\vect W_+(x,\cdot))(s) &= \mathcal{L}(u_-(0,\cdot))(s) \vect v_1(s)  e^{-\lambda(s)x} \,, (x,s)\in \Omega_+ \times \mathbb{C}^+
        \end{split}\right.
    \end{equation}
    with 
    \begin{equation}
        { \vect W}_0^- = (\eta_0^-, u_0^-) = \vect W_0 + H \mathcal{L}^{-1}(r) * G R \vect W_0\, \text{ in } \Omega
        \label{eq:hybrid-svb-solution-ic}
    \end{equation}
\end{remark}

\subsection{Well-posedness}
To simplify the notation let us write
\begin{equation}
    a \lesssim b
\end{equation}
to say that there exists a constant $C>0$ such that $a \leq C b$. Using \Cref{remark:hybrid-solution} we will now prove the well-posedness of the hybrid model for the case of data supported on one side of the interface. 

\begin{theorem}
    \label{th:continuity-one-side}
    Let $n \geq 0$, $\vect W_0 \in H^{n+1}(\Omega)\times H^{n+1}(\Omega)$ with $supp(\vect W_0) \subset \mathbb{R}^{-}$ and $\vect W$ the solution of the hybrid model. There exists a constant $C>0$ such that 
    \begin{equation}
        |\vect W|_{(H^{n,n}(Q))^2} \leq C |\vect W_0|_{H^{n+1}\times H^{n+1}(\Omega)}
    \end{equation}
\end{theorem}
\begin{proof}
    \emph{B to SV case:}
    Let $\mathcal{L}_- = \mathcal{L}_{B}$ and $\mathcal{L}_{+}=\mathcal{L}_{SV}$ and $\vect W_- = \vect W_{|Q_-}$ and $\vect W_+ = \vect W_{|Q_+}$. From \eqref{eq:exact-solution-abstract} we have $\vect W_- = \vect W^\star_{cauchy,-} + \vect W'_-$, so using \Cref{th:b-reflection} we deduce that $\vect W_-$ is solution of the Cauchy problem \eqref{eq:lgn-homogeneous-cauchy} with initial condition $\vect W_0^-$ from \eqref{eq:hybrid-bsv-solution-ic} and use that $|r|\leq 1$ and \Cref{th:cauchy-b}
    \begin{equation}
        \begin{split}
            |\vect W_-|_{H^{n,s}(Q_-)\times H^{n,s}(Q_-)} &\lesssim |\vect W_0^-|_{H^{n+1}(\Omega_-)\times H^{n+1}(\Omega_-)} \\
            & \leq |\vect W_0 |_{H^{n+1}(\Omega)\times H^{n+1}(\Omega)}
        \end{split}, \quad \forall s \geq 0
    \end{equation}
    and $u_-' \in H^{n+1,s}(Q_-)$ so from \Cref{th:trace} $u_-(0,\cdot) \in H^{\nu_0}$ with $\nu_0 = s \dfrac{n+1/2}{n+1}$. Choosing $s$ such that $\nu_0 = n$ then from \Cref{th:halfline-sv} and the continuity of the trace in \Cref{th:trace}
    \begin{equation}
        \begin{split}
            |\vect W_+|_{(H^{n,n}(Q_+))^2} &\lesssim |u_-(0,\cdot)|_{H^n} \\ 
            &\lesssim |\vect W_-|_{H^{n,n}(Q_-)\times H^{n,n}(Q_-)} \\ 
            &\lesssim |\vect W_0|_{H^{n+1}(\Omega)\times H^{n+1}(\Omega)}
        \end{split}
    \end{equation}
    hence reconstructing the solution we obtain 
    \begin{equation}
        |\vect W|_{(H^{n,n}(Q))^2}  \lesssim |\vect W_0|_{(H^{n+1}(\Omega))^2}
    \end{equation}

    \emph{SV to B case:}
    Let now $\mathcal{L}_- = \mathcal{L}_{SV}$ and $\mathcal{L}_{+}=\mathcal{L}_{B}$. From \Cref{th:sv-reflection}, $\vect W_-$ is the solution of the Cauchy problem \eqref{eq:sv-cauchy} with initial condition $\vect W_0^-$ given by \cref{eq:hybrid-svb-solution-ic}
    and use that $|r|\leq 1$ with  \Cref{th:sv-cauchy} to deduce
    \begin{equation}
        \begin{split}
            |\vect W_-|_{(H^{m,m}(Q_-))^2} &\lesssim |\vect W_0^-|_{H^{m}(\Omega_-)\times H^{m}(\Omega_-)} \\
            & \lesssim |\vect W_0 |_{H^{m}(\Omega)\times H^{m}(\Omega)}
        \end{split}
    \end{equation}
    also $u_-' \in H^{m,m}(Q_-)$ and $u_-(0,\cdot) \in H^{m}(0,T)$. Then, from \Cref{th:halfplane-boussinesq} and the continuity of the trace w.r.t the initial data (which can be derived from \eqref{eq:sv-cauchy-solution})
    \begin{equation}
        \begin{split}
            |\vect W_+|_{H^{r,m+1}(Q_+)\times H^{r,m}(Q_+)} &\lesssim |u_-(0,\cdot)|_{H^{m}} \\ 
            &\lesssim |\vect W_0^-|_{(H^{m}(\Omega_-))^2} \\ 
            &\lesssim |\vect W_0|_{(H^{m}(\Omega))^2} \\
        \end{split}
    \end{equation}
    Taking $r=m$ and reconstructing $\vect W|_{Q_-} = \vect W_-$ and $\vect W_{Q_+}=\vect W_+$ 
    \begin{equation}
        |\vect W|_{(H^{m,m}(\Omega))^2}  \lesssim |\vect W_0|_{H^{m}(\Omega)\times H^{m}(\Omega)}
    \end{equation}
    and taking $m=n+1$ this implies 
    \begin{equation}
        |\vect W|_{(H^{n,n}(Q))^2}  \lesssim |\vect W_0|_{(H^{n+1}(\Omega))^2}
    \end{equation}
\end{proof}
\begin{theorem}
    Let $n\geq 0$, if $\vect W_0 \in H^{n+1}(\Omega)\times H^{n+1}(\Omega)$ and $\vect W$ is the solution of the hybrid model, then there exists a constant $C>0$ such that 
    \begin{equation}
        |\vect W|_{(H^{n,n}(Q))^2} \leq C |\vect W_0|_{H^{n+1}\times H^{n+1}(\Omega)}
    \end{equation}
\end{theorem}
\begin{proof}
    The proof will be done by a density argument. Let $\vect W \in M\times M$, with 
    \begin{equation}
        M = \{f \in C^\infty(\mathbb{R}) : \quad f(0) = 0)\}
    \end{equation}
    and let $\mathcal{T}_{0}$ such that 
    \begin{equation}
        \vect W = \mathcal{T}_0 \vect W_0 
    \end{equation}
    is the solution of hybrid's problem \eqref{eq:hybrid-compact}. Define 
    \begin{equation}
        \vect W_0^L = \left\{\begin{array}{cc}
            \vect W_0 & \quad \text{ in } \Omega_- \\ 
             0 & \quad \text{ in } \Omega_+ 
        \end{array} \right.,\quad 
        \vect W_0^R = \left\{\begin{array}{cc}
             0  & \quad \text{ in } \Omega_- \\ 
            \vect W_0 & \quad \text{ in } \Omega_+ 
        \end{array} \right.,
        \label{eq:w0-splitting}
    \end{equation}
    then from \Cref{th:continuity-one-side}
    \begin{equation}
        |\mathcal{T}_0 \vect W_0^L |_{(H^{n,n}(Q))^2} \lesssim |\vect W_0^L|_{(H^{n+1}(\Omega))^2}
    \end{equation}
    \begin{equation}
        |\mathcal{T}_0 \vect W_0^R |_{(H^{n,n}(Q))^2} \lesssim |\vect W_0^R|_{(H^{n+1}(\Omega))^2}
    \end{equation}
    and by linearity
    \begin{equation}
        |\mathcal{T}_0 \vect W_0 |_{(H^{n,n}(Q))^2} \lesssim |\vect W_0|_{(H^{n+1}(\Omega))^2}
    \end{equation}
    which means that $\mathcal{T}_0$ is a linear bounded operator from $(C^\infty_0(\mathbb{R}))^2 \subset (H^{n+1}(\Omega))^2$ to $(H^{n,n}(Q))^2$. From \Cref{th:M0-dense-L2} we deduce that $(C^\infty_0(\mathbb{R}))^2$ is dense in $(H^{n+1}(\Omega))^2$ which, due to \Cref{th:unique-extension}, means that there exists a unique continuous operator 
    \begin{equation}
        \begin{array}{lll}
            \mathcal{T}:& (H^{n+1}(\Omega))^2 &\to (H^{n,n}(Q))^2\\
                & \vect W_0 & \to \vect W = \mathcal{T} \vect W_0
        \end{array}
    \end{equation}
    such that $\mathcal{T} \vect W_0 = \mathcal{T}_0 \vect W_0$ for every $\vect W_0 \in M$.    
\end{proof}

\subsection{Approximation and coupling error }
\label{sec:coupling-error}
Now we will discuss in an abstract sense the relation between the size of $\vect W'$, which we propose to name as the \emph{coupling error},  and the total approximation error $\vect W - \vect W_{3D}$. Then we will proceed to compute its asymptotic value for the case of the hybrid model.

Let  $\vect W^{\star\star}$ be the solution of the one-way coupling between the 3D model and the $\mathcal{L}_+$ model: 
\begin{equation}
    \vect W^{\star\star} = \left\{ 
    \begin{array}{cc}
        \vect W_{3D} & \text{ in } Q_- \\ 
        \vect W_{half,+}^{\star\star} & \text{ in } Q_+
    \end{array}
    \right.
    \label{eq:one-way-model-right-moving-3d}
\end{equation} 
with $\vect W_{half,+}^{\star\star}$ the solution of the half line problem on $Q_+$ 
\begin{equation} 
    \left\{ 
        \begin{array}{cc}
            \mathcal{L}_{+} \vect W_{+.half}^{\star\star}=0, & \quad \text{ in } Q_+ \\ 
            \vect W_{+.half}^{\star\star}(\cdot,0) = 0 & \quad \text{ in } \Omega_+ \\ 
            u_{half,+}(0,\cdot) = u^{3D}(0,\cdot) & \quad \text{ on } ]0,T[
        \end{array}
        \label{eq:one-way-3d-halfline-problem}
    \right., 
\end{equation}

Using these definitions we can separate the error in three parts as
\begin{equation}
    \vect W - \vect W^{3D} = (\vect W - \vect W^\star)  + (\vect W^\star - \vect W^{\star\star}) + (\vect W^{\star\star} -\vect W^{3D})
\end{equation}
which, in reverse order, can be interpreted as: 

\begin{itemize}
    \item \parbox[t]{\linewidth}{$\vect W^{\star\star} -\vect W^{3D}$ is  the \textbf{half-line model error} of the operator $\mathcal{L}_+$, because, from \eqref{eq:one-way-3d-halfline-problem}, one has 
    \begin{equation}
        \vect W^{\star\star} -\vect W^{3D} = \left\{\begin{array}{cc}
            0 & \text{ in } Q_- \\ 
            \vect W_{half,+}^{\star\star} & \text{ in } Q_+
        \end{array}\right.
    \end{equation}
    and since $\vect W_{half,+}^{\star\star}$ is only using data from the 3D model it does not carry errors from the other parts of the domain. It seems also that the most intuitive way to improve this error is by improving the operator $\mathcal{L}_+$, although the boundary condition \eqref{eq:one-way-3d-halfline-problem} could also influence the results. Such problems for BT equations have only recently been studied in the literature for arbitrary boundary data \cite{LannesWeynans2020,Bresch2021Waves} but they fall outside of the scope of the present study.
    }.
    
    \item \parbox[t]{\linewidth}{$\vect W^{\star} - \vect W^{\star\star}$ is the \textbf{Cauchy model error} of the $\mathcal{L}_-$ operator. To see this notice also that 
        \begin{equation}
            \vect W^{\star} - \vect W^{\star\star} = \left\{\begin{array}{cc}
                \vect W_{cauchy,-}^{\star} - \vect W_{3D} & \text{ in } Q_- \\ 
                \vect W_{half,+}^{\star} - \vect W_{half,+}^{\star\star} & \text{ in } Q_+
            \end{array}\right.
        \end{equation}
        so for $Q_-$ the comparison is w.r.t. to the best solution and for $Q_+$ the comparison is w.r.t. the same half-line model $\mathcal{L}_+$ but using the best possible data on the boundary, $u_{3D}(0,\cdot)$. In other words, in $Q_-$ this type of error indicates the precision the Cauchy problem of model $\mathcal{L}_-$, while for $Q_+$ it represents the part of this error that is carried by $\mathcal{L}_+$. Again, the most natural way of improving this error is by replacing $\mathcal{L}_-$ by a better model.
    \item \parbox[t]{\linewidth}{$\vect W - \vect W^\star$ is the \textbf{coupling error}. In this case we have 
        \begin{equation}
            \vect W - \vect W^{\star} = \left\{\begin{array}{cc}
                \vect W_- - \vect W_{cauchy,-}^{\star} & \text{ in } Q_- \\ 
                \vect W_+ - \vect W_{half,+}^{\star} & \text{ in } Q_+
            \end{array}\right.
        \end{equation}
        so for $Q_-$ this error measures the absorbing ability of the coupling conditions for the operator $\mathcal{L}_+$, while in $Q_+$ it describes how these errors propagate with the model $\mathcal{L}_+$. This coupling error will be analyzed in \Cref{sec:coupling-error}.
    }
    }
\end{itemize}
Even though the half-line model error and the Cauchy model error are also affected by the choice of transmission conditions, the one that more directly describes the effects added by the coupling method is the \textbf{coupling error}. The other types of error are also affected by transmission conditions, however, they seem more related to the specific problems on the whole and half lines.

\begin{figure}
    \centering
    \begin{subfigure}[b]{\textwidth}
        \centering
        \includegraphics[height=2in]{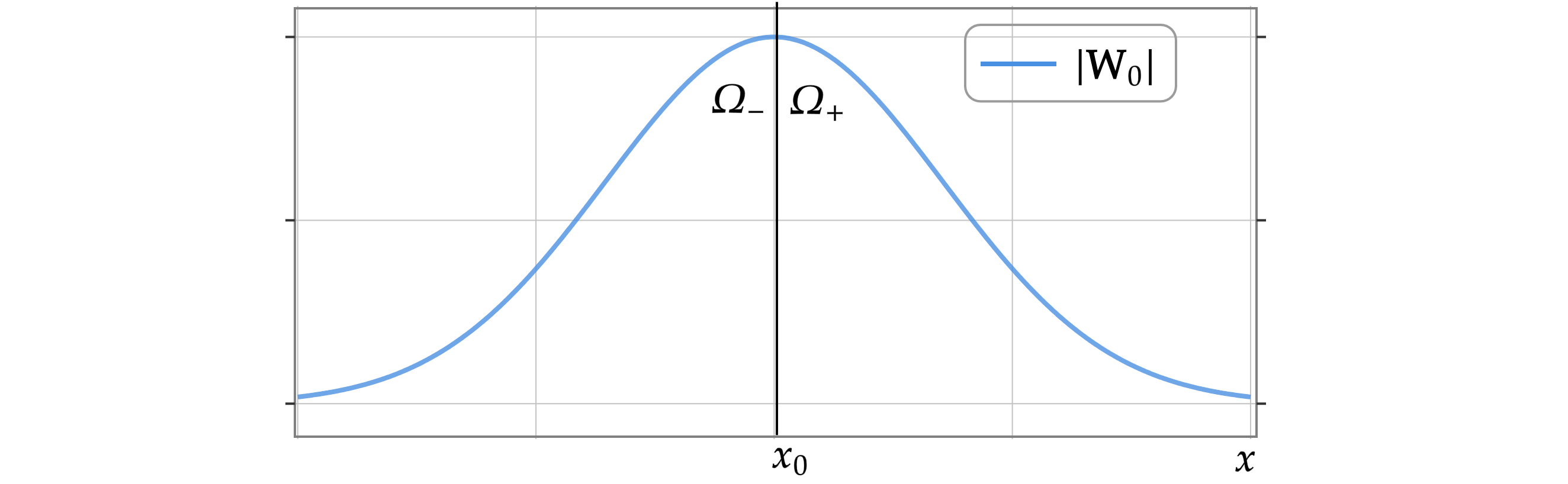}
    \end{subfigure}

    \vspace{1cm} %

    \begin{subfigure}[b]{\textwidth}
        \centering
        \includegraphics[height=2in]{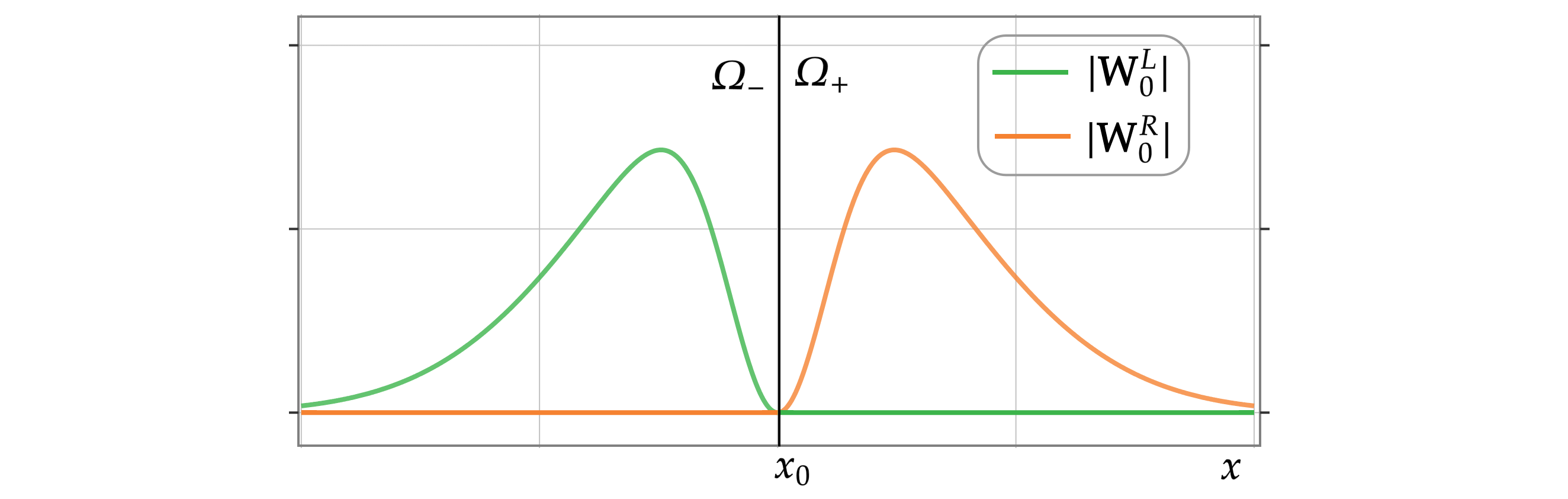}
    \end{subfigure}
    
    \caption{Sketch of the approximation used to handle the case of initial data with arbitrary support. The original initial condition $\vect W_0$ (top figure) is approximated as the superposition of two initial conditions $\vect W_0^L$, $\vect W_0^R$ that are supported on the left and on the right of the interface respectively.}
    \label{fig:w0w0lw0r}
\end{figure}

One important limitation of the one-way model is the assumption that the initial data is supported on one side of the interface. We will remove this assumption by approximating the initial condition $\vect W_0$ as the sum of two functions $\vect W_0^L$ and $\vect W_0^R$, each supported to the left and right of the interface, as shown in \Cref{fig:w0w0lw0r}. This approximation can be done arbitrarily thanks to a density argument that is proved in the \cref{app:approx}. The reference solution in this case is the superposition of two coupled models:

\begin{equation}
    \vect W^{\star\star\star}  = \vect W^\star_L  + \vect W^\star_R
\end{equation}
where $\vect W^\star_L$ is the solution of the one-way coupling \eqref{eq:one-way-model-right-moving} initialized with $\vect W_0^L$ as in \eqref{eq:w0-splitting}; and $\vect W^\star_R$ is the solution of the one-way coupling initialized with $\vect W_0^R$, but solved in the opposite direction, \textit{i.e.},
\begin{equation}
    \vect W^\star_R = \left\{ 
    \begin{array}{cc}
        \vect W_{cauchy,+} & \text{ in } Q_+ \\ 
        \vect W_{half,-} & \text{ in } Q_- 
    \end{array}
    \right.
\end{equation} 
where
\begin{equation}
    \left\{ 
        \begin{array}{cc}
            \mathcal{L}_{+} \vect W_{cauchy,+} =0, & \quad \text{ in } Q_+ \\ 
            \vect W_{cauchy,+} (\cdot, 0) = \vect W_0^R & \quad \text{ in } \Omega_+
        \end{array}
    \right.
\end{equation}
and
\begin{equation} 
    \left\{ 
        \begin{array}{cc}
            \mathcal{L}_{-} \vect W_{half,-} =0, & \quad \text{ in } Q_- \\ 
            \vect W_{half,-}(\cdot, 0) = 0 & \quad \text{ in } \Omega_+ \\ 
            u_{half,-}(0,\cdot) = u_{cauchy,+}(0,\cdot) & \quad \text{ on } ]0,T[
        \end{array}
    \right., 
\end{equation}

\subsubsection{Computation of the coupling error}
Now the goal is to analyze the coupling error $|\vect W'_\pm|$. The next Lemma \ref{th:trace-reflection-asymptotic-estimate} quantifies the asymptotic effect that the reflection coefficient has on a given function $f(t)$.

\begin{lemma}
    \label{th:trace-reflection-asymptotic-estimate}
    Let $f \in H^{n+2}(\mathbb{R}^+)$ with $n \geq 0$ integer, there exits $C>0$ such that
    \begin{equation}
        |\mathcal{L}^{-1}(r) * f |_{H^{n}(\mathbb{R}^+)} \leq C \mu^2 |f|_{H^{n+2}(\mathbb{R}^+)}
    \end{equation}
\end{lemma}
\begin{proof}
    From Theorem \ref{eq:approx-theory} there exists a family of functions $(\hat f_h)_{h>0} \subset C_0^\infty(\mathbb{R})$, and positive constants $c_1, c_2$, such that 
    $$\hat f_h = \hat f \text{ in } \left.\right] -\frac{c_1}{h},\frac{c_1}{h}\left[\right. \subset supp(\hat f_h) = \left.\right]  -\frac{c_2}{h}, \frac{c_2}{h}\left[\right.$$ and
    \begin{equation}
            |\hat f - \hat f_h|_{H^{n}(\mathbb{R})} \leq C_1 h |\hat f|_{H^{n+1}(\mathbb{R})}
            \label{eq:bound-for-fh}
    \end{equation}
    so we can compute 
    \begin{equation}
        \begin{split}
            |r \hat f|_{H^{n}(\mathbb{R})}^2 &\leq |r(\hat f-\hat f_h)|^2_{H^n(\mathbb{R})} + |r \hat f_h|^2_{H^n(\mathbb{R})}\\
            &= \int_\mathbb{R} (1+\omega^2)^{n/2}|r|^2 |\hat f -\hat f_h|^2 d\omega + \int_\mathbb{R} (1+\omega^2)^{n/2} |r|^2 |\hat f_h|^2 d\omega \\ 
            & =\int_{|\omega|\geq \frac{c_1}{h}} (1+\omega^2)^{n/2}|r|^2 |\hat f - \hat f_h|^2 d\omega  + \int_{|\omega| \leq \frac{c_2}{h}} (1+\omega^2)^{n/2}|r|^2 |\hat f_h|^2 d\omega \\ 
        \end{split}
    \end{equation}
    where the limit of the first integral is obtained from the region where $\hat f = \hat f_h$, and the limit of the second integral from the support of $\hat f_h$. Since $|r| \leq 1$, we can proceed to use \eqref{eq:bound-for-fh} to obtain 
    \begin{equation}
        |r \hat f|_{H^n(\mathbb{R})}^2 \leq C_1^2 h^2 |\hat f|_{H^{n+1}(\mathbb{R})}^2 + \int_{|\omega| \leq \frac{c}{h}} |r|^2 |\hat f_h|^2 d\omega
    \end{equation}
    Since we have
    \begin{equation}
        r(s) = -\frac{1}{4}\mu^2 s^2 + O(|\mu s|^4) \quad \text{ if } |\mu s| \leq  c' < 1
    \end{equation}
    then choosing $h = \dfrac{c}{c'} \mu^2$ and $\sigma$ very small, developing the integral one obtains that there exists $C_2>0$ such that 
    \begin{equation}
        \int_{|\omega| \leq \frac{c}{h}} |r|^2 |\hat f_h|^2 d\omega \leq C_2 \mu^4 |f_h|_{H^{n+2}}
    \end{equation}
    from where we deduce that there exists $C_3$ such that 
    
    \begin{equation}
        |\mathcal{L}^{-1} (r) * f|_{H^n(\mathbb{R}^+)}^2 \leq \mu^4  C_3 |f_h|_{H^{n+2}(\mathbb{R}^+)}^2
    \end{equation}
\end{proof}

The next \Cref{th:asymptotic-size-of-reflections-onesided} shows that as $\mu\to 0$ the size of the coupling error  $\vect W_\pm'$, \textit{i.e.}, the reflections, is $O(\mu^2)$.

\begin{proposition}
    \label{th:asymptotic-size-of-reflections-onesided}
    Let $\vect W_0 \in (H^n(\Omega))^2$ and $supp(\vect W_0)\subset \Omega_-$, for every $r, k > 0 $ there exists $C(r,k)>0$ and $C>0$ such that  if $\mathcal{L}_- =\mathcal{L}_{B}$ then 
    \begin{equation}
        \begin{split}
            |\vect W_-'|_{H^{r,k+1}(Q_-)\times H^{r,k}(Q_-)} \leq C(r,k) \mu^2 |\vect W_0|_{(H^n(\Omega))^2} \\ 
            |\vect W_+'|_{H^{k,k}(Q_+)\times H^{k,k}(Q_+)} \leq C(r,k) \mu^2 |\vect W_0|_{(H^n(\Omega))^2} \\ 
        \end{split}
    \end{equation}
    and if $\mathcal{L}_-  = \mathcal{L}_{SV}$ then there is $C_1>0$ and $C_2(r)>0$ such that 
    \begin{equation}
        \begin{split}
            |\vect W_-'|_{H^{n,n}(Q_-)\times H^{n,n}(Q_-)} \leq C_1 \mu^2 |\vect W_0|_{(H^{n+2}(\Omega))^2} \\ 
            |\vect W_+'|_{H^{r,n+1}(Q_+)\times H^{r,n}(Q_+)} \leq C_2(r) \mu^2 |\vect W_0|_{(H^{n+2}(\Omega))^2} \\ 
        \end{split}
    \end{equation}
\end{proposition}
\begin{proof}
    Suppose $\mathcal{L}_- = \mathcal{L}_B$ then for all $r \geq 0$ and $k\geq 0$
    \begin{equation}
        \begin{split}
            |\vect W_-'|_{H^{r,k+1}(Q_-)\times H^{r,k}(Q_-)} &\lesssim |u_-'(0,\cdot)|_{H^k}(0,T)  \quad \text{(\Cref{th:halfplane-boussinesq})}\\ 
            &\lesssim C \mu^2 |u^\star_{cauchy,-}(0,\cdot)|_{H^{k+2}(0,T)} \quad \text{(\Cref{th:trace-reflection-asymptotic-estimate})} \\ 
            &\lesssim C \mu^2 |\vect W_{cauchy,-}^\star|_{(H^{n,s_k}(Q))^2} \quad \text{(\Cref{th:trace})} \\ 
            &\lesssim \mu^2 | \vect W_0|_{(H^n(\Omega))^2} \quad \text{(\Cref{th:cauchy-b})}
        \end{split}
    \end{equation}
    with $s_k = (k+2)\frac{n+1}{n+1/2}$, thanks to the continuity of the trace in \Cref{th:trace}. Similarly
    \begin{equation}
        \begin{split}
            |\vect W_+'|_{H^{k,k}(Q_+)\times H^{k,k}(Q_+)} &\lesssim |u_-'(0,\cdot)|_{H^k}(0,T)  \quad \text{(\Cref{th:halfline-sv})}\\ 
            &\lesssim C \mu^2 |u^\star_{cauchy,-}(0,\cdot)|_{H^{k+2}((0,T))} \quad \text{(\Cref{th:trace-reflection-asymptotic-estimate})} \\ 
        &\lesssim C \mu^2 |\vect W_{cauchy,-}^\star|_{(H^{n,s_k}(Q))} \quad \text{(\Cref{th:trace})} \\ 
            &\lesssim \mu^2 | \vect W_0|_{(H^n(\Omega))^2} \quad \text{(\Cref{th:cauchy-b})}
        \end{split}
    \end{equation}
    Suppose now that $\mathcal{L}_- = \mathcal{L}_{B}$, then 
    \begin{equation}
        \begin{split}
            |\vect W_-'|_{H^{n,n}(Q_-)\times H^{n,n}(Q_-)} &\lesssim |u_-'(0,\cdot)|_{H^n}(0,T)  \quad (\text{\Cref{th:halfline-sv}})\\ 
            &\lesssim C \mu^2 |u^\star_{cauchy,-}(0,\cdot)|_{H^{n+2}(0,T)} \quad \text{(\Cref{th:trace-reflection-asymptotic-estimate})} \\ 
            &\lesssim C \mu^2 |\vect W_{cauchy,-}^\star|_{(H^{n+2,n+2}(Q))^2} \quad \text{(\Cref{th:sv-cauchy})} \\ 
            &\lesssim \mu^2 | \vect W_0|_{(H^{n+2}(\Omega))^2} \quad \text{(\Cref{th:cauchy-b})}
        \end{split}
    \end{equation}
    and for every $r \geq 0$
    \begin{equation}
        \begin{split}
            |\vect W_+'|_{H^{r,n+1}(Q_+)\times H^{r,n}(Q_+)} &\lesssim |u_-'(0,\cdot)|_{H^n(0,T)} \quad \text{ (\Cref{th:halfplane-boussinesq})}\\ 
            &\lesssim C \mu^2 |u^\star_{cauchy,-}(0,\cdot)|_{H^{n+2}(0,T)} \quad \text{(\Cref{th:trace-reflection-asymptotic-estimate})} \\ 
            &\lesssim C \mu^2 |\vect W_{cauchy,-}^\star|_{(H^{n+2,n+2}(Q))^2} \quad \text{(\Cref{th:trace})} \\ 
            &\lesssim \mu^2 | \vect W_0|_{(H^{n+2}(\Omega))^2} \quad \text{(\Cref{th:cauchy-b})}
        \end{split}
    \end{equation}    
\end{proof}
To generalize this result to the case of initial data with arbitrary support, let us define 
\begin{equation}
    \vect W^{\star\star\star}  = \vect W^\star_L  + \vect W^\star_R
\end{equation}
where $\vect W^\star_L$ is the solution of the one-way coupling \eqref{eq:one-way-model-right-moving} initialized with $\vect W_0^L$ as in \eqref{eq:w0-splitting}; and $\vect W^\star_R$ is the solution of the one-way coupling initialized with $\vect W_0^R$, but solved in the opposite direction, \textit{i.e.},
\begin{equation}
    \vect W^\star_R = \left\{ 
    \begin{array}{cc}
        \vect W_{half,-} & \text{ in } Q_- \\
        \vect W_{cauchy,+} & \text{ in } Q_+ \\ 
    \end{array}
    \right.
\end{equation} 
where
\begin{equation}
    \left\{ 
        \begin{array}{cc}
            \mathcal{L}_{+} \vect W_{cauchy,+} =0, & \quad \text{ in } Q_+ \\ 
            \vect W_{cauchy,+} (\cdot, 0) = \vect W_0^R & \quad \text{ in } \Omega_+
        \end{array}
    \right.
\end{equation}
and
\begin{equation} 
    \left\{ 
        \begin{array}{cc}
            \mathcal{L}_{-} \vect W_{half,-} =0, & \quad \text{ in } Q_- \\ 
            \vect W_{half,-}(\cdot, 0) = 0 & \quad \text{ in } \Omega_+ \\ 
            u_{half,-}(0,\cdot) = u_{cauchy,+}(0,\cdot) & \quad \text{ on } ]0,T[
        \end{array}
    \right., 
\end{equation}
so defining $\vect W' = \vect W - \vect W^{\star\star\star}$ one has the following result:

\begin{theorem}
    \label{th:asymptotic-size-of-reflections}
    Let $\vect W_0 \in H^{n+2}(\Omega) \times H^{n+2}(\Omega)$, there exists $C>0$ such that as $\mu$ goes to 0
    \begin{equation}
        |\vect W'|_{(H^{n,n}(\Omega))^2} \leq C \mu^2|\vect W_0|_{H^{n+2}(\Omega)\times H^{n+2}(\Omega)}
    \end{equation}
\end{theorem}
\begin{proof}
    Let 
    \begin{equation}
        M = \{f \in C^\infty(\mathbb{R}) : \quad f(0) = 0)\}
    \end{equation}
    and suppose $\vect W_0 \in M$. Also notice that 
     \begin{equation}
        \begin{split}
            \vect W' &=  \vect W - \vect W^{\star\star\star}  \\
            &= (\vect W_L - \vect W_L^\star) + (\vect W_R - \vect W_R^\star)
        \end{split}
    \end{equation}
    and from \Cref{th:asymptotic-size-of-reflections-onesided} one can deduce that 
    \begin{equation}
        |\vect W_L'|_{(H^{n,n}(Q))^2} \leq C \mu^2 |\vect W_0|_{H^{n+2}(\Omega)\times H^{n+2}(\Omega)}
    \end{equation}
    \begin{equation}
        |\vect W_R'|_{(H^{n,n}(Q))^2} \leq C \mu^2 |\vect W_0|_{H^{n+2}(\Omega)\times H^{n+2}(\Omega)}
    \end{equation}
    which means that 
    \begin{equation}
        |\vect W'|_{(H^{n,n}(Q))^2} \leq C \mu^2 |\vect W_0|_{H^{n+2}(\Omega)\times H^{n+2}(\Omega)}
        \label{eq:asymptotic-behavior}
    \end{equation}
    which means that the mapping 
    \begin{equation}
        \begin{array}{ll}
            \mathcal{T}'_0: M^2 &\to (H^{n,n}(Q))^2  \\ 
            \vect W_0 &\to \vect W'
        \end{array}
    \end{equation}
    is continuous, and because $M$ is dense in $H^{n+2}(\Omega)$ (\Cref{th:M0-dense-L2}) it has a unique extension  (\Cref{th:unique-extension})
    \begin{equation}
        \begin{array}{ll}
            \mathcal{T}': (H^{n+2}(\Omega))^2 &\to (H^{n,n}(Q))^2  \\ 
            \vect W_0 &\to \vect W'
        \end{array}
    \end{equation}
    and \eqref{eq:asymptotic-behavior} holds.
\end{proof}

\begin{remark}[Directionality]
Despite the fact just proved that $\vect W' =O(\mu^2)$ in both cases (SVB and BSV), one must notice in \eqref{th:asymptotic-size-of-reflections-onesided} that when the wave propagates in the BSV case the regularity estimate of $\vect W'$ is much better than in the SVB case. Moreover, after examining the filter $r(j\omega)$ in both \eqref{th:b-reflection} and \eqref{th:sv-reflection}, in the limit $\sigma \to 0$, these  filters take the form 
\begin{equation}
r_{BSV}(j\omega_1(\kappa)) = \dfrac{\sqrt{1+(\mu\kappa)^2}-1}{\sqrt{1+(\mu\kappa)^2}+1}    
\label{eq:bsv-filter}
\end{equation}
and 
\begin{equation}
    r_{SVB}(j\omega_1(\kappa))= \dfrac{1-\sqrt{1-(\mu\kappa)^2}}{1+\sqrt{1-(\mu\kappa)^2}}    
    \label{eq:wbt-filter}
\end{equation}
one has that $|r_{SVB}(\kappa)|$ and $|r_{BSV}(\kappa)|$ are increasing functions of $\mu|\kappa|$ but $|r_{SVB}(\kappa)|> |r_{BSV}(\kappa)|$. Morover, in the limit $ |\mu \kappa|\to \infty$ $|r_{SVB}| =|r_{BSV}| =1$ but $|r_{BSV}|<1$ for every $\kappa$, while $|r_{SVB}|=1$ for every $\frac{h_0^2}{3} \kappa^2>1$. This means that one should expect larger reflections in deeper waters (larger $\mu$) but even larger for the SVB case.
\end{remark}

%% file: A01halfline.tex
\label{sec:appendix}
Here we proceed to analyze problem \eqref{eq:half-line-problem}, proving the well-posedness  \Cref{th:halfplane-boussinesq} and extending the results of references \cite{Fokas2005,johnston2021,mantzavinos2023} by removing the assumptions of smooth data and compatibility conditions to data in $H^s(0,T)$ and no compatibility conditions.  %

\subsection{Well-posedness of the Boussinesq equations on the halfline}

First we need some lemmas:
\begin{lemma}
    \label{th:sqrt-newton}
    Let $f(z) = (1+z)^{-1/2}$ and $|z|\leq r < 1$. Then $f(z) = 1+R(z)$ with $|R(z)| \leq C |z|$ for some  $C>0$.
\end{lemma}
\begin{proof}
First, let $z=\rho e^{i \theta}$, $\rho >0$ and $\theta \in [0, 2\pi)$, then $|f(z)|^2 \leq \dfrac{1}{1-\rho} = f(-\rho)$, so by Taylor expansion of $f(z)$, letting $M=f(-\rho)$, for every $|z|<r<1$
    $$
        |R(z)| \leq \dfrac{M|z|}{\rho(\rho-|z|)}\leq \frac{M|z|}{\rho(\rho-r)}= C |z| \\
    $$
\end{proof}
\begin{lemma}
    Let $\mu>0$ and $0<r<1$, there exists $C>0$ and $R(s)$ such that for every $s \in \mathbb{C}^+$ and $|s|^2 > \dfrac{1}{\mu^2 r}$:
    $$
        \lambda(s) = \frac{1}{\mu} + R(s), \quad |R(s)|  \leq C \frac{1}{\mu^2|s|^2}, \quad \text{ and } Re(R(s)) > 0
    $$
    and $C$ only depends on $r$.
    \label{th:lambda-large-s}
\end{lemma}
\begin{proof}
Since we can write
    $$ \lambda(s) = \dfrac{1}{\mu} \left(1+\dfrac{1}{\mu^2 s^2}\right)^{-1/2} $$ 
then because of the assumption $|1/(\mu s)|<r<1$ we can apply Lemma \ref{th:sqrt-newton} and obtain such $C$ and $R(s)$. Moreover, because $0<Re(\lambda(s)) < |\lambda(s)| < \dfrac{1}{\mu}$, we also have $Re(R(s))>0$.
\end{proof}
\begin{lemma}
    Let $\mu>0$, there exists $\sigma_0>0$ such that if $|s|> \sigma_0$ then
    $$
        \frac{1}{2\mu} < |\lambda(s)| < \frac{3}{2\mu}
    $$
    and
    $$
        \frac{1}{2\mu} < Re(\lambda(s)) < \frac{3}{2\mu}
    $$
    \label{th:eigenvalues-bounds}
\end{lemma}
\begin{proof}
Let $C$ and $R(s)$ be as in \Cref{th:lambda-large-s}, then if $|s|> \sigma_0 = 2C /\mu$ and $|R(s)| \leq \frac{1}{2\mu}$, the conclusion is obtained from $$\dfrac{1}{\mu} - |R(s)| \leq |\lambda(s)| \leq \dfrac{1}{\mu} + |R(s)|$$ and $$\dfrac{1}{\mu} - |R(s)| \leq Re(\lambda(s)) \leq \dfrac{1}{\mu} + |R(s)|$$.
\end{proof}

Now we use the previous lemmas for Theorem \eqref{th:halfplane-boussinesq}.

\begin{proof}[Proof of Theorem \eqref{th:halfplane-boussinesq}]
    The proof is done by a density argument.  Let $$M=\{ f \in C_0^\infty(]0,T[): \int_{\mathbb{R}} f(\omega) d\omega = 0\}$$ and $M_0$ the space of functions $f :[0,T] \to \mathbb{R}$ that can be written as
    \begin{equation}
        f(t) = e^{\sigma t} \int_\mathbb{R} \alpha_\sigma( \omega) e^{j\omega t} d\omega, \quad \alpha_\sigma \in M
        \label{eq:f-alpha-M0}
    \end{equation}
    with $\sigma > 0$. One can see that because of the assumptions on $\alpha_\sigma$ one has $f \in C^\infty(]0,T[)$, $f(0)=0$ and $\alpha_\sigma(\omega)$ is the Laplace transform of $f$ at $s=\sigma+j\omega$. We will prove that for any $n \geq 0$, $\vect W_-$ given by \eqref{eq:half-line-solution} defines a linear bounded operator 
    \begin{equation}
        \left\{
        \begin{array}{rl}
            M_0 \mapsto& H^{n,k+1}(Q_-)\times H^{n,k}(Q_-) \\ 
            \mathcal{T}_0 u_\Gamma =& \vect W_-
        \end{array}\right.
    \end{equation} 
    and then use \Cref{th:unique-extension} to show that there is a unique extension to
    \begin{equation}
        \left\{
        \begin{array}{rl}
             H^{k}(0,T) \mapsto& H^{n,k+1}(Q_-)\times H^{n,k}(Q_-) \\ 
            \mathcal{T} u_\Gamma =& \vect W_-
        \end{array}\right.
        \label{eq:halfline-extended-operator}
    \end{equation} 

    Let $\sigma>0$ fixed and $f,\alpha$ as in \eqref{eq:f-alpha-M0}. Letting $u_{\Gamma-}=f$ in \eqref{eq:half-line-solution} 
    \begin{equation}
        \widehat{\vect W}_-(x, \sigma+j\omega) = \alpha_\sigma(\omega) \vect v_2(\sigma+j\omega) e^{\lambda(\sigma+j\omega)} x
    \end{equation}
    so taking the inverse Laplace transform 
    \begin{equation}
        \vect W_-(x, t) = e^{\sigma t}\int_\mathbb{R}\alpha_\sigma(\omega) \vect v_2(\sigma+j\omega) e^{\lambda(\sigma+j\omega)x} d\omega
    \end{equation}
    we will verify that $\vect W_-$  verifies the initial and boundary conditions, proving that there exists a classical solution to \eqref{eq:half-line-problem}. 

    A direct evaluation at $x=0$ confirms the boundary condition $\vect W_-(0,\cdot) = f$. For the initial condition, because $\alpha_\sigma \in C_0^\infty$ and $\sigma>0$ we have $\vect W_-(x,\cdot) \in C^\infty(]0,T[)$ so the integral can be evaluated pointwise in $t$. Evaluating at $t=0$ we have 
    \begin{equation}
        \vect W_-(x,0) = \int_\mathbb{R} \alpha_\sigma(\omega) \vect v_2(\sigma+j\omega) e^{\lambda(\sigma+j\omega) x} d\omega
    \end{equation}
    so using Lemma \ref{th:lambda-large-s}
    \begin{equation}
        \begin{split}
            \vect W_- (x,0) &= e^{x/\mu} \int_\mathbb{R} \alpha_\sigma(\omega) \vect v_2 e^{R(s)x}d\omega \\
            &\approx e^{x/\mu}\int_\mathbb{R} \alpha_\sigma(\omega) \vect v_2 (1+R(s)x) d\omega \\ 
            &= e^{x/\mu}\int_\mathbb{R} \alpha_\sigma(\omega) \vect v_2  R(s)x dx
        \end{split}
    \end{equation}
    where we used the approximation of the exponential in the second line, and the assumption of the integral of $\alpha_\sigma$ in the third one. Now, using Lemma \ref{th:eigenvalues-bounds} again we can compute
    \begin{equation}
        |\vect W_-(x,0)| \leq e^{x/\mu} \int_\mathbb{R} |R(s) x||\vect v_2| d\omega \leq e^{x/\mu} \frac{2C |x|}{\mu^2\sigma^2}
    \end{equation}
    because $|\vect v_2| \leq 2$,  which means that for fixed $x \in \Omega_-$, since $\sigma>0$ can be arbitrarily large,  then necessarily $\vect W_-(x,0)=0$. This proves that indeed $\vect W_-$ is a solution of \eqref{eq:half-line-problem} in $Q_-$. Moreover, because the Laplace-transform representation \eqref{eq:half-line-solution} is unique, it is the unique solution.

    \vspace{1pt}
    Now, to prove the regularity, let $s = \sigma + j \omega$, because of \eqref{th:eigenvalues-bounds} we have that  if $\sigma>\sigma_0$ then $|s|>\sigma_0$ and

  \begin{equation}
      \dfrac{1}{2\mu} < | Re(\lambda(s))| < \dfrac{3}{2\mu}
  \end{equation}
  for every $\omega \in \mathbb{R}$. Now because of \eqref{eq:half-line-solution} we have that
  \begin{equation}
      |\eta_-(s)|  \leq 
      \frac{|\alpha_\sigma(s)|}{|\mu s||\sqrt{1+(\mu s)^{-2}}|}
      e^{Re(\lambda(s))x}, \quad |u_-(s)|  
      \leq 
      |\alpha_\sigma(s)|  e^{Re(\lambda(s))x}
  \end{equation}
  but there exists a lower bound $0<l(\mu, \sigma) \leq  |\sqrt{1+(\mu s)^{-2}}|$ so also
  \begin{equation}
    |\eta_-(s)|  \leq 
    \frac{|\alpha_\sigma(s)|}{ l(\mu, \sigma) \mu |s|}
    e^{Re(\lambda(s))x}
  \end{equation}
  We will first show a continuity estimate for $u$, but because of the previous inequalities the same steps also apply for $\eta_-$ and $\partial_t \eta_-$. From  \Cref{th:eigenvalues-bounds} we know there is $\sigma_0>0$ such that $\sigma > \sigma_0$ implies
  \begin{equation}
      Re(\lambda(s)) \geq \frac{1}{2\mu}
  \end{equation}
  so there is $C_1$ such that 
  \begin{equation}
      |{\widehat u_-}(s,x) | \leq C_1 |\alpha_-(s)| e^{\frac{1}{2\mu}x}, \quad x \in \Omega_-
  \end{equation}
  so for $k\geq 0$, integrating we obtain
  \begin{equation}
      \label{eq:trace-uniform-bound}
      \begin{split}
          \int_0^T \int_{\Omega^-} |{ \partial_t^k u_-} |^2 dx dt
          &\leq C(T) \int_{\Omega_-} \int_\mathbb{R} |s^k\hat u_-(\sigma + j\omega, x)|^2 d\omega dx \\
          &\leq C(T) C_1 \int_\mathbb{R} |s^k \alpha_-(\sigma + j\omega)|^2 \int_{-\infty}^0 e^{x/\mu} dx d\omega \\
          &= C(T) C_1 \mu \int_\mathbb{R} |s^k \alpha_-(\sigma + j\omega)|^2 d\omega\\
          &= C(T)C_1 \mu \int_0^T |f^{(k)}(t)|^2 dt
      \end{split}
  \end{equation}
  where the last line comes from the causality principle. This shows that $\eta_-, u_-$ and $\partial_t \eta_-$ are in $H^{0,k}(Q_-)$ and 
  \begin{equation}
      |\vect W_-|_{H^{0,k+1}(Q_-)\times H^{0,k}(Q_-)} \leq C(T)C_1\mu |f|_{H^k(0,T)}
  \end{equation}

  Also, for $n$ positive and integer $$\partial_x^n \vect{\widehat W_-} = \lambda(s)^n \vect{\widehat W_-}$$ so, using Theorem \eqref{th:eigenvalues-bounds}
  $$
      |\partial_x^n \vect{\widehat W_-}| = |\lambda(s)|^n |\vect{\widehat W_-}| \leq \left(\frac{3}{2\mu}\right)^n |\vect{\widehat W_-}|
  $$
  so by Parseval's theorem  $\partial_x^n \vect {W _-}\in H^{0,
  k+1}\times H^{0,k}$, \textit{i.e.}, $\vect W_- \in H^{n,k}$, and
  \begin{equation}
      | \partial_x^n \vect W_-|_{L^2(0,T;L^2(\Omega_-))} \leq C(T)C_1\mu \left(\frac{3}{2\mu}\right)^n|f|_{L^2(0,T)}
  \end{equation}
  from where we obtain $C_2(n,\mu)>0$ such that 
  \begin{equation}
    | \vect W_-|_{H^{n,k+1}(Q_-)\times H^{n,k}(Q_-)} \leq C(T)C_1\mu C_2(n,\mu) |f|_{H^k(0,T)}
  \end{equation}
  So far, we have proved that the operator 
  \begin{equation}
    \left\{
        \begin{array}{ll}
             M_0 \subset H^{k}(0,T) &\to H^{n,k+1}(Q_-)\times H^{n,k}(Q_-) \\ 
            \mathcal{T}_0 f &= \vect W_-
        \end{array}\right.
  \end{equation}
  is well-defined, linear and bounded. Because of Theorem \ref{th:M-dense-L2} we know that $M_0$ is a dense subspace of $H^k(]0,T[)$ so by Theorem \ref{th:unique-extension} this means that there exists a unique bounded operator $\mathcal{T}: H^{k}(0,T) \to H^{n,k+1}(Q_-)\times H^{n,k}(Q_-)$ such that $ \mathcal{T} f = \mathcal{T}_0 f $  for any $f \in M_0$.

  Now, writing $ \vect{\tilde W}_-=(\tilde \eta, \tilde u) =( \partial_t^k\eta, \partial_t^k u)$, then from \eqref{eq:half-line-problem}
  \begin{equation}
      \begin{split}
          \partial_t \tilde \eta + \partial_x \tilde u = 0 \\
          (1-\mu^2\partial_x^2)\partial_t \tilde u +  \partial_x \tilde \eta = 0
      \end{split}
  \end{equation}
  which can be rewritten as
  \begin{equation}
      \begin{split}
          (1-\mu^2 \partial_x^2)\partial_t^2 \tilde u  -  \partial_x^2 \tilde u = 0 \\
          (1-\mu^2 \partial_x^2)\partial_t^2 \tilde \eta  -  \partial_x^2 \tilde \eta = 0
      \end{split}
  \end{equation}
  by taking $\partial_x$ of the first one, $\partial_t$ of the second and adding them; and by  taking $(1-\mu^2\partial_x^2)\partial_t$ of the first one and $(1-\mu^2 \partial_x^2)\partial_x$ of the second one and adding them. So by  \Cref{th:elliptic} we obtain that $\partial_t^2 \tilde u \in H^{n,k}(Q_-)$ and $\partial_t^2 \tilde \eta \in H^{n,k+1}(Q_-)$ and, since $T<\infty$, we deduce that $\vect {\tilde W} \in H^{k,2}$, hence $\vect W \in H^{n,k+3}(Q_-)\times H^{n,k+1}(Q_-)$ for any nonnegative integer $k$.
\end{proof}